\documentclass{siamltex}
\usepackage{graphicx}
\usepackage{amsmath}
\usepackage{amssymb}

\def\Real{\mathbb R}
\def\R2{{\mathbb R}^2}

\def\ua{u_\alpha}

\author{
David C. Dobson \and
Lyubima B. Simeonova\thanks{
Both authors: Department of Mathematics, University of Utah,
Salt Lake City, UT 84112-0090, USA}
}

\title{Optimization of periodic composite structures for sub-wavelength
       focusing\thanks{
This work was partially supported by NSF grant DMS-0537015.}
}

\begin{document}
\maketitle

\begin{abstract}
Recently, there has been plenty of work in designing and fabricating
materials with an effective negative refractive index. Veselago realized
that a slab of material with a refractive index of $-1$ would act as a lens.
Pendry suggested that the Veselago lens would act as a {\em superlens},
providing
a perfect image of an object in contrast to conventional lenses which are only
able to focus a point source to an image having a diameter of the order of
the wavelength of the incident field. 
 
Recent work has shown that similar focusing effects can be obtained
with certain slabs of ``conventional'' periodic composite materials: 
photonic crystals.
The present work seeks to answer the question of what periodic dielectric
composite medium (described by dielectric coefficient with positive 
real part) gives an optimal image of a point source. An optimization 
problem is formulated and
it is shown that a solution exists provided the medium
has small absorption.  Solutions are characterized by an adjoint-state 
gradient condition, and several numerical examples illustrate
both the plausibility of this design approach, and the possibility of
obtaining smaller image spot sizes than with typical photonic crystals.
\end{abstract}

\pagestyle{myheadings}
\thispagestyle{plain}
\markboth{D. C. DOBSON AND L. B. SIMEONOVA}
{OPTIMIZATION OF SUB-WAVELENGTH FOCUSING}

\section{Introduction}

Recently, there has been a renewed and avid interest in studying a class of materials known as the
left-handed materials (LHMs). These materials have simultaneously negative real parts of dielectric permittivity $\rho$ and 
magnetic permeability $\mu$, so that their refractive index is negative. The properties of such meterials were investigated first by Veselago in 1967 \cite{Veselago}. 
As shown by Veselago, LHMs exhibit some peculiar electromagnetic properties such as negative index of
refraction and  wave vector, $\textbf{k}$, and Poynting vector, $\textbf{S}$, having opposite
directions. Veselago realized that a slab of LHM would act as a lens. 

%Rays emanating from a source would be 
%focused first inside the slab and then refocused on the opposite side of the lens. 

Due to the absence of naturally ocurring materials possessing both negative permittivity and 
negative permeability, Veselago's predictions did not receive much attention until recently,
when a material with both negative permittivity and negative permeability at microwave frequencies was built \cite{Smith}. 
Subsequently, the properties of LHMs were analyzed by many authors. 

According to Abbe's diffraction limit, conventional lenses based on positive index materials with curved surfaces
are not able to resolve an object's fine details that are smaller than half of the light wavelength $\lambda$. 
The limitation occurs because the waves with transverse wave numbers larger than $2 \pi n/\lambda$, which carry 
information about the fine sub-$\lambda$ details of the object, decay exponentially in free space. 
In a negative index material slab, however, the evanescent wave components can grow exponentially 
and thus compensate for the exponential decay. Therefore, under ideal conditions, all Fourier components from
the object can be recovered at the image plane producing a resolution far below the diffraction limit \cite{Pendrytwo}.

%However, the far field "perfect" lens may only have limitted applications because any realistic losses can eliminate the 
%superlensing effect \cite{Podolskiy}. Moreover, the achievable resolution is limited by the absorption and thickness
%of the LHM slab \cite{Rao}. 
%Kildishev \textit{et al.} apply stochastic optimization tools to design a low-loss optical 
%negative index metamaterial \cite{Kildishev}. 

Milton \textit{et al.} proved superlensing in the quasistatic
regime (where the wavelength is much larger than the object), and discused limitations of superlenses in this 
regime due to anomalous localized resonance. If the source being imaged responds to an applied field, it must lie
outside the resonant regions to be successfully imaged \cite{Milton}. 

In the electrostatic limit, the magnetic and electric fields decouple, and the requirement for superlensing of 
transverse-magnetic waves is reduced to only $\rho = - \rho_h$, where $\rho_h$ is the permittivity of the host medium 
interfacing the lens \cite{Pendrytwo}. An example of such near field superlens is a slab of silver in air illuminated 
at its surface plasmon resonance (where $\rho = -1$). Experiments with silver slabs have already shown rapid growth 
of evanescent waves \cite{Liu}, submicron imaging \cite{Melville}, and imaging beyond the diffraction limit \cite{Fang}.
A major draw-back of such near-field superlenses based on bulk metals is that they can operate only at a single frequency 
$\omega$ satisfying the lens condition $\rho (\omega) = - \rho_h$. Shalaev \textit{et al.} proposed a "tunable" 
near-field superlens made of metal-dielectric composites that can operate at any desired visible or near-infrared wavelength
with the frequency controlled by the metal filling factor of the composite (here the inhomogeneities are assumed to be much 
smaller than the wavelength) \cite{Shalaev}.
 
%The concept of superlens was questioned by some authors
%and arose some controversy (\cite{Maystre}, \cite{Garcia}). Garcia and Nieto-Vesperinas showed that the field generated by a point source 
%through such a material is described by an integral that diverges in some regions of space. 
%They concluded that this integral "cannot represent any physically realizable wave field," and claimed that sufficient 
%losses inside the material "drastically change any evanescent amplifying wave into a decaying one" (\cite{Garcia}).

It was shown by Efros \textit{et al.} that a two dimensional photonic crystal made from a non-magnetic dielectric
has negative values of both the electric permittivity and the magnetic permeability in some frequency range
\cite{Efros_one}, and the photonic crystal behaves like a real LHM with respect to the propagating modes only.
 If amplification of evanescents waves occurs, it is due to some other reason; for example, excitation of surface waves by
the evanescent waves. Such an amplification may provide an improvement of the image in the near-field region, but
it does not affect the image near the far-field focal point \cite{Efros_two}.

The physical principles that allow negative refraction in photonic crystals arise from the dispersion
characteristics of wave propagation in periodic media and are very different from those in LHMs. They also do not
require both negative electric permittivity and magnetic permeability \cite{Kosaka,Notomi}.
The negative refraction of beams can be described by analyzing the equifrequency surface of the band
structures \cite{Kosaka,Notomi,Luo}. If the constant-frequency contour is everywhere convex, an
incoming plane wave from air will couple to a single mode that propagates into the crystal on the negative side of
the boundary, and thus negative refraction in the first band is realized.  Luo \textit{et al.} have shown all-angle 
negative refraction could be achieved at the lowest band of two-dimensional photonic crystals in the case of 
$\textbf{S} \cdot\textbf{k} > 0$ \cite{Luo}. Such all-angle refraction is essential for superlens application. 
The photonic crystal not only focuses all propagating waves without limitation of finite aperture, but also
amplifies at least some evanescent waves, and the unconventional imaging effects are due to the
presence of additional near-field light. A perfect lens, made of left-handed materials, focuses all
propagating waves and all evanescent waves. The important difference for superlensing with a
photonic crystals is that only finite number of evanescent waves is amplified. This is a consequence of Bragg
scattering of light to leaky photon modes \cite{Luotwo}. The resolution of a photonic-crystal superlens at a single
frequency is only limited by its surface period instead of the wavelength
\cite{Luotwo}.

%All-angle negative refraction and evanescent wave
%amplification  is achieved also by Shin ans Fan using a one-dimensional photonic crystal consisting
%of metal-dielectric multilayers (\cite{Shin}). 

More recently Huang \textit{et al.} proposed an
alternative approach to all-angle negative refraction in two-dimensional photonic crystals. By
applying appropriate modifications with surface grating to the flat photonic lens, he is able to focus large
and/or far way objects \cite{Huang}.

Inspired by the current research in structures that produce sub-wavelength focusing, we use derivative-based
minimization techniques to produce structures that will provide sub-wavelength focus with non-magnetic materials and without the need for negative 
permittivities. 
Rather than restricting to designs based on photonic crystal structures,
we allow as admissible {\em any} periodic composite structure (with
fixed period) whose refractive index is bounded above and below by
fixed constants.  And rather than performing parametric optimization over
a small number of variables describing the structure, we use techniques
of ``topology optimization'' in which material distribution is completely
arbitrary. We are able to obtain structured that focus a point source
that is far away from the lens, and also we can obtain structures that give an
image at a chosen distance from the lens.  Since structures incorporating
gratings are included in our admissible class, such designs 
will naturally arise through the optimization process if they
produce the best possible image.

For simplicity, only the case of ``two-dimensional'' structures
in $E$-parallel polarization is considered.  The ideas here should
extend to the other polarization case
and the full three-dimensional problem, although there are some
technical hurdles. 

The paper proceeds as follows. In Section 2 we describe the model problem and review a variational formulation 
of the Helmholtz equation in a periodic geometry.  The inclusion of a 
small amount of energy absorption in the medium allows a uniform
upper bound on the norm of the electric field, independent of the
particular admissible structure (and thereby preventing resonances).
 In Section 3 we present the optimization problem and derive the 
optimality conditions. In Section 4 we depict the numerical experiments and show structures that have produced sub-wavelength 
focus. 

\section{Model problem}
In this paper we consider time-harmonic electromagnetic wave propagation through nonmagnetic $(\mu = 1)$
heterogeneous media for which the dielectric coefficient is constant in one direction, i.e.
$\epsilon(x,y,z)= \rho(x,y)$. Assuming that the electric field vector $E=(0,0,u)$, Maxwell's equations 
reduce to the Helmholtz equation
\begin{equation} \label{helmholtz}
\triangle u + \omega^2 \rho u = 0, \quad \mbox{in $\R2$},
\end{equation}
where $\omega$ represents the frequency, and
$\rho \in L^\infty(\R2)$
is the dielectric coefficient.

\subsection{Periodic structure}
Assume that the dielectric coefficient $\rho(x,y)$ is
periodic in the $x$ variable
\[ \rho (x,y) = \rho (x+2\pi,y), \quad 
                             \mbox{for all $(x,y) \in \R2$}.
\]
Taking the period to be $2 \pi$ imposes no loss of generality since
any other period can be obtained by rescaling $\omega$.

Assume that the regions $\{ y > 0 \}$, and
$\{ y < -b \}$ are  homogeneous, for some fixed constant $b > 0$.
In particular, assume for $y > 0$ and $y < -b$, that
$\rho(x,y) = 1$.  The slab $-b < y < 0$ may
contain inhomogeneous material.

Suppose a point source is placed above the slab at the point
$(0,h)$, which generates the incident field 
$u_i(x,y) = H_0^{(1)}(\omega r)$, where 
$r = \sqrt{x^2 + (h - y)^2}$ is the distance from the source,
and $H_0^{(1)}$ is the Hankel function.  For $y < h$, we have
the representation 
\[ u_i(x,y) = \frac{1}{\pi} \int_\Real \frac{1}{\beta(\xi)}
e^{i \xi x - i \beta(\xi) (y - h)} \, d\xi,
\]
where $\beta(\xi) = \sqrt{\omega^2 - \xi^2}$ whenever
the argument is positive, 
and $\beta(\xi) = i \sqrt{\xi^2 - \omega^2}$ otherwise \cite{Morse}.
It follows that
\[ f(x) \equiv u_i(x,0) = 
\frac{1}{\pi} \int_\Real \frac{e^{i\beta(\xi)h}}{\beta(\xi)}
e^{i \xi x} \, d\xi,
\]
and
\[ g(x) \equiv \frac{\partial u_i}{\partial y}(x,0) = 
\frac{-i}{\pi} \int_\Real e^{i\beta(\xi)h} e^{i \xi x} \, d\xi.
\]
We can rewrite
\[ f(x) = \int_{-1/2}^{1/2} f_\alpha(x) e^{i\alpha x} \, d\alpha,
\quad \mbox{where } \;\;
f_\alpha(x) = \frac{1}{\pi} \sum_{n \in {\mathbb Z}}
\frac{e^{i\beta(n+\alpha )h}}{\beta(n + \alpha )} e^{inx},
\]
and
\[ g(x) = \int_{-1/2}^{1/2} g_\alpha(x) e^{i\alpha x} \, d\alpha,
\quad \mbox{where } \;\;
g_\alpha(x) = \frac{-i}{\pi} \sum_{n \in {\mathbb Z}}
   {e^{i\beta(n+\alpha )h}} e^{inx}.
\]
($f_\alpha$ may fail to converge, but only at the isolated
values of $\alpha$ for which $(n+\alpha)^2 = \omega^2$ for some $n$.
$g_\alpha$ always converges due to the exponential decay in $n$.)

Above the slab $y > 0$, we separate the solution $u$ 
to (\ref{helmholtz})
into the incident and scattered field: $u = u_i + u_s$.
The scattered field $u_s$ can also be separated by the Fourier
transform in $x$,
\[ u_s(x,y) = \int_\Real \hat{u}_s(\xi,y) e^{i\xi x}\, d\xi.
\]
Plugging this representation back into (\ref{helmholtz}),
and solving for each frequency $\xi$ separately, we find
that 
$\hat{u}(\xi,y) = a(\xi)e^{i\beta(\xi) y} + b(\xi) e^{-i\beta(\xi) y}$.
The second term on the right corresponds to an {\em incoming} wave,
which we insist must be zero since we want the scattered field to
consist only of outgoing waves.  Then
$\hat{u}_s(\xi,y) = a(\xi) e^{i\beta(\xi)y}$.
It follows that $u_s(x,0) = \int a(\xi) e^{i\xi x} \, d\xi$,
and 
\begin{eqnarray*}  \frac{\partial u_s}{\partial y}(x,0) &=& 
 \int i\beta(\xi) a(\xi) e^{i\xi x} \, d\xi \\
 &=& \int i\beta(\xi) \hat{u}_s(\xi,0) e^{i\xi x} \, d\xi
 \equiv (Tu_s)(x).
\end{eqnarray*}
The linear operator $T$ (Dirichlet-to-Neumann map) 
then defines the relationship between the traces
$u_s|_{\{y=0\}}$ and $\partial_y u_s|_{\{y = 0\}}$:
$T(u_s|_{\{y = 0\}}) = (\partial_y u_s)|_{\{y=0\}}$.
On the boundary $\{ y = 0\}$, the solution
$u = u_i + u_s$ should then satisfy
\begin{eqnarray*}
 \partial_y u - Tu &=& \partial_y u_i - Tu_i + \partial_y u_s -Tu_s
 = g -Tf  \\
 &=& 2g.
\end{eqnarray*}

Define the periodic domain (circle) 
\[ \mathcal S = \Real/2 \pi \mathbb Z. \]
Define the first Brillouin zone $K=[-\frac{1}{2}, \frac{1}{2}]$. To reduce the problem (\ref{helmholtz}) over $\R2$
to a family of problems over $\mathcal S \times \Real$, we define for $g \in L^2(\R2)$ the Floquet transform
$\mathcal F$ by
\[\mathcal F (g) = e^{-i\alpha x} \sum_{n \in Z} g(x-2 \pi n,y)e^{i 2 \pi \alpha n}, \quad \alpha \in K. \]
The sum can be considered as a Fourier series in the quasi-momentum variable $\alpha$, with values in 
$L^2(\mathcal S \times \Real)$. The map $g \mapsto \mathcal F g$ is an isomorphism from $L^2(\R2)$ to the direct
product space $\int^\oplus_K L^2(\mathcal S \times \Real)$ \cite{Kuchment}.
Floquet theory assures that the solution
$u$ can be written
\begin{equation} \label{uexp}
 u(x,y) = \int_{-1/2}^{1/2} u_\alpha(x,y)e^{i \alpha x} \, d\alpha,
\end{equation}
where each function $u_\alpha$ is $2\pi$-periodic in the $x$ 
variable, and satisfies the equation
\[ \triangle_\alpha u_\alpha + \omega^2 \rho u_\alpha = 0,
\]
where $\triangle_\alpha = \triangle + 2 i \alpha \partial_1 - |\alpha|^2$.
The boundary condition 
$\partial_y u - T u = 2g$ on $\{y = 0\}$ translates to
\[
 \partial_y \ua  - T_\alpha \ua = 2 g_\alpha, \quad \mbox{on $\{y = 0\}$},
\]
where
\begin{equation} \label{Tdef}
T_\alpha \ua = \sum_{n\in {\mathbb Z}} 
 i \beta(n + \alpha) \hat{u}_\alpha(n) e^{inx},
\end{equation}
(here $\hat{u}_\alpha(n)$ are the Fourier series coefficients of $\ua$ with respect to the $x$-variable).
Similar considerations apply at the lower boundary of the slab
$\{ y = -b \}$, where (assuming there is no incoming 
wave coming from below) we find 
$\partial_y \ua + T_\alpha \ua = 0$.

\subsection{Existence and uniqueness of solutions}
Let $\Omega = \mathcal S \times (0, -b)$, $\Gamma_0 = \{y = 0 \}$, $\Gamma_b = \{y = -b \}$.
Define an admissible class of dielectric coefficients 
\[ \mathcal{A} = \{ \rho = \rho_r + i\rho_{i} \in L^\infty(\Omega): 
  \rho_{r_0} \leq \rho_r(x) \leq \rho_{r_1}
  \mbox{ and }\rho_{{i}_0} \leq \rho_{i}(x) \leq \rho_{{i}_1} \; \mbox{a.e.}\},
\] 
where $\rho_{{r}_0}$, $\rho_{{i}_0} > 0$.
Given the incident wave $u_i$ generated by the point source
at $(0,h)$, we must solve the family of problems
\begin{align} \label{feqn}
\triangle_\alpha \ua + \omega^2 \rho_r \ua + i\omega^2 \rho_{i} \ua & =  0,  
                                       &\mbox{in $\Omega$}\\
\nonumber
(\frac{\partial}{\partial y} - T_\alpha) \ua & = 2g_\alpha,
                        &\mbox{ on $\Gamma_0$} \\
\nonumber
(\frac{\partial}{\partial y} + T_\alpha) \ua & = 0,
                        &\mbox{ on $\Gamma_b $},
\end{align}
for all $\alpha \in [-\frac{1}{2}, \frac{1}{2}]$.
Existence and uniqueness of weak solutions, with a uniform bound, may be
obtained for $\rho_{i} > 0$.

\begin{lemma} \label{Solution}
For each $\rho \in  {\mathcal A}$ with $\rho_{i} > 0$ and $\alpha \in [-\frac{1}{2}, \frac{1}{2}]$, problem  (\ref{feqn})
admits a unique weak solution $u_\alpha \in H^2(\Omega)$.  Furthermore,
there exists a constant $C$ depending on $\rho_{i}$, ${\mathcal A}$,
such that $\| u_\alpha \|_{H^2(\Omega)} \leq C$, 
independent of $\rho \in  {\mathcal A}$ and $\alpha$.
\end{lemma}
\begin{proof}
For convenience we drop the subscript $\alpha$ on solutions.
Define for $u$, $v \in H^1(\Omega)$
\[ a(u,v) = \int_\Omega \nabla u \cdot 
               \overline{\nabla v} - 
               \omega^2\int_\Omega \rho u \overline{v} + \alpha^2\int_\Omega  u \overline{v} - 2i\alpha \int_\Omega \partial_x {u} \overline{v} -
	        \int_{{\Gamma}_0} (T_{\alpha} u) \overline{v} - \int_{{\Gamma}_b} (T_{\alpha} u) \overline{v}, 
\]
and
\[ b(v) =   2\int_{{\Gamma}_0} g_{\alpha} \overline{v} .
\]
It is straightforward to show
that $a(u,v)$ defines a bounded sesquilinear form over
$H^1(\Omega) \times H^1(\Omega)$, and that
$b(v)$ is a bounded linear functional on $H^1(\Omega)$.
Weak solutions $u \in H^1(\Omega)$ of (\ref{feqn}) solve the
variational problem
\begin{equation} \label{varprob}
a(u,v) = b(v) \quad \mbox{for all $v \in H^1(\Omega)$}.
\end{equation}
The sesquilinear form $a$ uniquely defines a linear operator $A: H^1(\Omega) 
\rightarrow H^1(\Omega)$
such that $a(u,v) = \langle Au, v \rangle_{H^1(\Omega)}$, 
and the functional $b(v)$ is uniquely
identified with an  element $b \in H^1(\Omega)$ such 
that $b(v) = \langle b, v \rangle$
by reflexivity and an abuse of notation.  Problem (\ref{varprob}) is
then equivalently stated
\begin{equation} \label{abst_varprob}
A u = b.
\end{equation}

We intend to show that $a$ is coercive by establishing a bound 
$|a(u,u)| \geq c > 0$
for all $u \in H^1(\Omega)$ with $\| u \|_{H^1(\Omega)} = 1$.

Integrating by parts in $x$, we see by periodicity that 
 \[ \int_\Omega (\partial_x {u}) \overline{u} = 0 \]
 
 Let $\Lambda^+ (\alpha) = \lbrace n \in {\mathbb Z} \colon Im (\beta_n) = 0 \rbrace$ and $\Lambda^-(\alpha) = {\mathbb Z} - \Lambda^+ (\alpha)$.
 Each $\Lambda^+ (\alpha)$ is a finite set and $0 \in \Lambda^+ (\alpha)$. 
  We can write
 \[ - \int_{{\Gamma}_j} (T_{\alpha} u) \overline{u} = 
 - \int_{\Gamma_j} \sum_{n \in \Lambda^+ (\alpha)} i \beta_n(\alpha)\hat{u}(n) e^{i nx} \overline{u} 
 - \int_{\Gamma_j} \sum_{n \in \Lambda^- (\alpha)} i \beta_n(\alpha)\hat{u}(n)
 e^{ inx} \overline{u},  \]
for $j = 0$, $b$. Notice that all $\beta_n(\alpha) \in \Lambda^- (\alpha)$ satisfy $- i \beta_n \geq 0$, 
so the second term on the right hand side of the equation above is real and non-negative.
  We then have
\begin{eqnarray*} 
a(u,u)   =  \int_\Omega |\nabla u|^2 + \alpha^2 \int_\Omega |u|^2 +
\int_{\Gamma_j} \sum_{n \in \Lambda^- (\alpha)} -i\beta_n(\alpha)\hat{u}(n) e^{inx} \overline{u} 
 - \omega^2\int_\Omega \rho_r |u|^2\\ 
 - i\int_{\Gamma_j} \sum_{n \in \Lambda^+ (\alpha)} \beta_n(\alpha)\hat{u}(n) e^{i nx} \overline{u} 
 -i\omega^2\rho_{i} \int_\Omega |u|^2.
    \end{eqnarray*}
Assuming $\|u\|_{H^1(\Omega)}^2 = 
    \int_\Omega |\nabla u|^2 +\int_\Omega |u|^2  = 1$, 
and noticing that
the first four terms on the right-hand side are purely 
real and the last two terms
are purely imaginary, we find
\begin{eqnarray*}
2| a(u,u) |   \geq  \left| 1 
+ \int_{\Gamma_j} \sum_{n \in \Lambda^- (\alpha)} -i\beta_n(\alpha)\hat{u}(n) e^{inx} \overline{u} 
                   - \int_\Omega (1 + \omega^2 \rho_r - {\alpha}^2 ) |u|^2 \right| \\
         + \left| \int_{\Gamma_j} \sum_{n \in \Lambda^+ (\alpha)} \beta_n(\alpha)\hat{u}(n) e^{i nx} \overline{u} + 
\omega^2 \rho_{i} \int_\Omega |u|^2 \right|.
\end{eqnarray*}
For convenience, write 
$t = \int_{\Gamma_j} \sum_{n \in \Lambda^- (\alpha)} -i\beta_n(\alpha)\hat{u}(n) e^{inx} \overline{u} $,
$r = \int_\Omega(1 + \omega^2 \rho_r \}) |u|^2$, and
$s = \int_\Omega |u|^2$.
Obviously $t,r$, and $s$ are nonnegative real numbers which depend on $u$ (and
$\rho$ in the case of $r$). 
Although $t$ and $s$ are essentially independent, $r$ must satisfy
\begin{equation} \label{rbounds}
(1 + \rho_{r_0})s \leq r \leq (1 + \rho_{r_1}) s.
\end{equation} 
With this notation,
\[ 2|a(u,u)| \geq |1 + t +{\alpha}^2 s - r| + \omega^2 \rho_{i} s.
\]
Note that in the case $s \geq  \frac{1}{2(1+\rho_{r_1})}$,  we have 
$|a(u,u)| \geq \frac{1}{2}\omega^2 \rho_{i} s \geq \frac{\omega^2 \rho_{i}}{4(1+\rho_{r_1})}$.
Otherwise, $s <  \frac{1}{2(1+\rho_{r_1})}$ so that $r < \frac{1}{2}$, and
$|a(u,u)| \geq \frac{1}{2}|1 + t - r|  > \frac{1}{4}$.
Hence, for all $s, t \geq 0$, and all $r$ satisfying (\ref{rbounds}),
\[ |a(u,u)| \geq  c = 
\min\left\{\frac{\omega^2 \rho_{i}}{4(1+\rho_{r_1})}, \frac{1}{4} \right\}.
\]
The bound thus holds for every $u$ with $\| u \|_{H^1(\Omega)} = 1$ 
and for every 
$\rho \in {\mathcal A}$ with $\rho_{i} > 0$.
Given this coercivity bound, direct application of the Lax-Milgram 
Theorem yields
existence of the bounded  solution operator 
$A^{-1}$ for problem (\ref{abst_varprob})
such that $\| A^{-1} \| \leq 1/c$.  
Thus $\| u \|_{H^1(\Omega)} \leq \| b \|_{H^1(\Omega)}/c$.

Given the  bound on $\| u \|_{H^1(\Omega)}$, a uniform $H^2(\Omega)$ bound
follows easily, since $\triangle_{\alpha} u = - \omega^2 \rho u$ 
is uniformly bounded in $L^2(\Omega)$.
\end{proof}

The complete solution $u$
to the original problem (\ref{helmholtz}) can then be
reconstructed from (\ref{uexp}).

\section{Optimal design}

The goal of the optimization is to make a perfect image of
the incident field $u_i=H_0^{(1)}(\omega r)$ on the opposite side of the slab
$\{ y < -b \}$.  A ``mirror image'' to the incident field 
converging at the point $(0, -(b+h_1))$ would look like
$H_0^{(2)}(\omega\sqrt{x^2 + (y+b+h_1)^2})$, where 
$H_0^{(2)} = \overline{H_0^{(1)}}$
is the conjugate Hankel function.  Thus we want the trace
\[ u(x,-b) = H_0^{(2)}(\omega\sqrt{x^2 + h_1^2}) \equiv q(x).\]

\begin{figure}
\centerline{
\includegraphics[width=70mm]{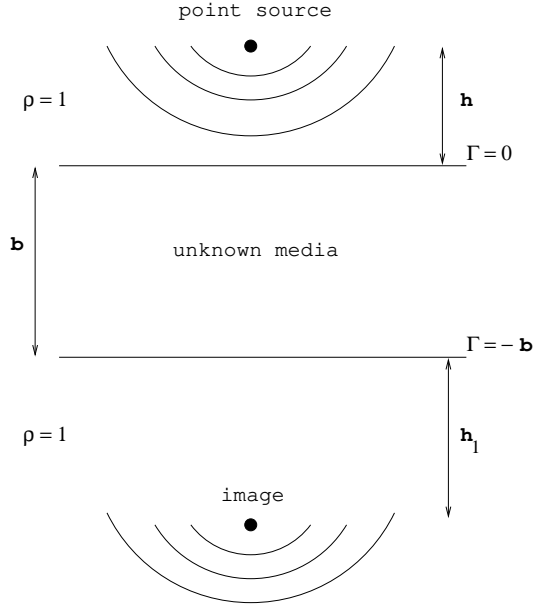}
}
\label{figure1}
\caption{Model problem.  A time-harmonic wave from a point source
is incident from above.  We wish to determine the unknown periodic
medium such that a focus is attained below.}
\end{figure}

The Bloch representations of $u$ and $f$ allows us to see that by setting
\begin{equation} \label{d_cond}
 \ua(x,-b) = q_\alpha(x),
\end{equation}
with $q_\alpha$ defined similarly to $f_\alpha$, we get $u(x,-b) = q(x)$.

Problem (\ref{feqn}) together with the additional boundary
condition (\ref{d_cond}) is overposed.  However, by allowing
$\rho$ to vary as a design variable, it may be possible to
make (\ref{d_cond}) hold approximately for each $\alpha$.

\subsection{Problem Formulation}

Let $F(\rho,\alpha) = \ua|_{\Gamma_b}$, 
where $\ua \in H^1(\Omega)$ is the
weak solution to problem (\ref{feqn}).

Consider the minimization
\begin{equation} \label{minprob}
 \inf_{\rho \in \mathcal{A}} J(\rho) = 
           \frac{1}{2} \int_{-1/2}^{1/2}
\| F(\rho,\alpha) - q_\alpha \|^2_2 \, d\alpha.
\end{equation}

\begin{theorem} \label{Optimization}
The optimization problem has a solution.
\end{theorem}
\begin{proof}
The proof follows the well-known direct method in the calculus of varations.
$\mathcal A$ is weak $\star$ $\ L^{\infty}$ compact. Consider a minimizing sequence $\lbrace \rho_n \rbrace$ with some subsequence 
(still denoted by $\lbrace \rho_n \rbrace$) converging weak $\star$ $\ L^{\infty}$ to some $\rho \in \mathcal A$.
Denote by $u_n$ the solution to the boundary value problem corresponding to $\rho_n$. By Lemma (\ref{Solution}), 
the sequence $\lbrace u_n \rbrace$ has has a weakly convergent
subsequence (still denoted $\lbrace u_n \rbrace$), $u_n \rightharpoonup u$ in $H^1(\Omega)$ for some $u \in H^1(\Omega)$. 
Hence, $u_n \longrightarrow u$ strongly in $L^2(\Omega)$. We hold $v \in H^1(\Omega)$ fixed, 
and we have $u_n\overline{v} \longrightarrow u \overline{v}$ strongly in $L^1(\Omega)$.
For $k$ fixed, $a_{\rho_k} (u_n,v) \longrightarrow a_{\rho_k} (u,v)$ as $n \longrightarrow \infty$. 
Here we used the fact that $T^{\alpha}:H^{1/2}(\Gamma_j) \to H^{-1/ 2}(\Gamma_j)$ is continuous. 
Since $u \overline{v} \in L^1(\Omega)$, $a_{\rho_k} (u,v) \longrightarrow a_{\rho} (u,v)$ 
as $\rho_k \longrightarrow \rho$ weak $\star$ $L^\infty$. 
The trace map extends uniquely to a continuous linear operator $\tau:H^1(\Omega) \longrightarrow H^{1/2}(\Gamma)$. 
Thus, the traces are also convergent: 
$u_n|_{\Gamma_j} \rightharpoonup u|_{\Gamma_j}$ weakly in $H^{1/2}(\Gamma)$. 
This implies $F(\rho_n, \alpha) \rightharpoonup F(\rho, \alpha)$ weakly in $H^{1 /2}(\Gamma)$, 
and, hence, $F_\alpha: \mathcal A \longrightarrow H^{1/2}(\Gamma)$ is weak $\star$ $L^\infty$ continuous for fixed $\alpha$.
But our bound on the solution is independent of $\alpha$, and thus, the minimization problem has at least one solution $\rho \in \mathcal A$.  
\end{proof}

\subsection{Adjoint-state derivatives}

Let $\delta \rho = \delta \rho_r + i \delta \rho_{i}$ be a "small" perturbation to the coefficient 
$\rho$. We denote the
linearization of $J(\rho)$ with respect to $ \delta \rho$ by $DJ(\rho)\delta \rho$. We have
\begin{eqnarray*} 
 DJ(\rho) =& \int_{-\frac{1}{2}}^{\frac{1}{2}} Re \langle DF_{\alpha}(\rho)\delta \rho, F(\rho,\alpha) -
q_{\alpha} \rangle_{L^2(\Gamma_b)} \, d\alpha \\
=& \int_{-\frac{1}{2}}^{\frac{1}{2}} Re \langle \delta \rho, DF_{\alpha}^{\star}(\rho) ( F(a,\alpha) -
 q_{\alpha} )\rangle_{L^2(\Omega)} \, d\alpha.
\end{eqnarray*}
For fixed $\alpha$, \[ DF_{\alpha}(\rho): L^2(\Omega)
\longrightarrow L^2(\Gamma_b) \]

and
\[DF_{\alpha}^{\star}(\rho): L^2(\Gamma_b)
\longrightarrow L^2(\Omega).\]

$DF_{\alpha}(\rho)(\delta \rho) = \delta u_{\alpha}\mid_{\Gamma_b}$ where $\delta u_{\alpha}$ solves the linearized problem:
\begin{align} \label{perturbprob}
\triangle_\alpha \delta u_{\alpha} + \rho \delta u_{\alpha} & =  -\delta \rho u,  
                                       & \mbox{in $\Omega$}\\
\nonumber
(\frac{\partial}{\partial y} - T_\alpha) \delta\ua & = 0,
                        & \mbox{on $\Gamma_0$} \\
\nonumber
(\frac{\partial}{\partial y} + T_\alpha) \delta\ua & = 0,
                        & \mbox{on $\Gamma_b$},
\end{align}

The $L^2$ adjoint of the derivative $DF_{\alpha}(\rho)(\cdot)$ is the linear operator $DF_{\alpha}^{\star}(\cdot)$ such that 
\[\langle DF_{\alpha}(\rho)(\delta \rho), \psi \rangle_{L^2(\Gamma_b)} = \langle \delta \rho, DF^{\star}(\rho)(\psi)\rangle_{L^2(\Omega)}.\]
For all $\psi \in L^2(\Gamma_b)$ let $w_{\alpha}$ solve

\begin{align} \label{wprob}
\triangle_\alpha w_{\alpha} + \omega^2\rho_r w_{\alpha} - i\omega^2 \rho_{i} w_{\alpha} & = 0,  
                                       & \mbox{in $\Omega$}\\
\nonumber
(\frac{\partial}{\partial y} - T^{\star}_\alpha) w_{\alpha} & = 0,
                        & \mbox{on $\Gamma_0$} \\
\nonumber
(\frac{\partial}{\partial y} + T^{\star}_\alpha) w_{\alpha} & = -\psi,
                        & \mbox{on $\Gamma_b$},
\end{align}
where $T^{\star}_\alpha f = -\sum i \overline{\beta}_n \hat{f}(n) e^{inx}$.
An integration by parts argument shows that 
\[ \int_{\Gamma_b} \delta u_{\alpha} = \int_{\Omega} \delta \rho u_{\alpha}
\overline {w}_{\alpha}.\]
The $L^2$ "gradient" of the functional $J(\rho)$ is the function $G(\rho) \in L^2(\Omega)$ for which
\[ DJ(\rho)(\delta \rho) = Re \left ( \int_{\Omega} \delta \rho G(\rho) \right ) \quad \mbox{and} \quad
G(\rho) = \int_{-\frac{1}{2}}^{\frac{1}{2}}u_{\alpha}\overline{w}_{\alpha} \, d\alpha \]
where $w_{\alpha}$ solves (\ref{wprob}) with $\psi = F (\rho,\alpha) - q_\alpha$.

 Define the normal cones at the minimizer $\hat {\rho}= \hat{\rho_r} + i \hat{\rho_{i}}$:
\[ N(\hat{\rho_r}) = \left \lbrace \xi \in (L^\infty)': \int_\Omega (\hat{\rho_r} - \rho_r)d \xi \geq 0 \quad
\forall \rho_r \in \mathcal A \right \rbrace,\]
and
\[ N(\hat{\rho_{i}} )= \left \lbrace \zeta \in (L^\infty)': 
\int_\Omega (\hat{\rho_{i}} - \rho_{i}) d \zeta \geq 0 \quad
\forall \rho_{i} \in \mathcal A \right \rbrace.\]
The gradient is normal to $\mathcal A$ at $\hat{\rho}$: $-G \bigcap N(\hat{\rho}) \not= 0$.
Thus, we have
\[ - \int_\Omega (\hat{\rho_r} - \rho_r) Re ( G(\rho)) dx \geq 0 \quad \forall \rho_r \in \mathcal A. \]
The continuity of the gradient is sufficient to reduce the above to pointwise optimality conditions. In
particular,
\begin{align*}
\hat{\rho_r} = \rho_{r_{0}} \quad  &\Rightarrow \quad Re(G(\rho)) > 0 \\
\rho_{r_{0}} < \hat{\rho_r} < \rho_{r_{1}} \quad  &\Rightarrow \quad Re(G(\rho)) = 0 \\
\hat{\rho_r} = \rho_{r_{1}} \quad  &\Rightarrow \quad Re(G(\rho)) < 0 
\end{align*}
for almost every x in $\Omega$.

Similarly, 
\[\int_\Omega (\hat{\rho_{i}} - \rho_{i}) Im ( G(\rho)) dx \geq 0 \quad \forall \rho_{i} 
\in \mathcal A,\]
and the pointwise optimality conditions are:
\begin{align*}
\hat{\rho_{i}} = \rho_{i_{0}} \quad  &\Rightarrow \quad Im(G(\rho)) < 0 \\
\rho_{i_{0}} < \hat{\rho_{i}} < \rho_{i_{1}} \quad  &\Rightarrow \quad Im(G(\rho)) = 0 \\
\hat{\rho_{i}} = \rho_{i_{1}} \quad  &\Rightarrow \quad Im(G(\rho)) > 0 
\end{align*}
for almost every x in $\Omega$.

In a similar optimal design problem involving waveguides, one can 
use the weak continuation property to show
\cite{Dobson} that at an optimal $\rho$,
at least one of the four equality bound constraints 
above on the real and imaginary parts of $\rho$ is attained at almost
every point $x \in \Omega$.  A similar argument for the present problem
is greatly complicated by the form of $\nabla J$, involving an integral
of a family of PDE solutions, rather than just a single background-adjoint
pair.  A precise a-priori characterization of optimal solutions is thus
difficult. 

%Suppose 
%\[ \lbrace x \colon \rho_{r_{0}} < \hat{\rho_r} < \rho_{r_{1}} \mbox{ and }
% \rho_{i_{0}} < \hat{\rho_{i}} < \rho_{i_{1}} \rbrace \]
% contains an open set $U$. \newline
% \centerline{}
% Recall that 
% $\nabla_\rho J(\rho) = \int_{-1/2}^{1/2}u_{\alpha}\overline{w}_{\alpha} \, d\alpha $, 
% where $u_\alpha$ and $w_\alpha$ are solutions of the direct and adjoint scattering problems.
% Since $u_\alpha, w_\alpha \in H^2(\Omega)$ are continuous with continuous derivatives, 
% one of the two must also be zero on  an open set contained in $U$. 
% By weak unique continuation property of solutions, 
% the field must then be zero everywhere. This can only occur when 
% the boundary data is zero, which is only possible for $w_\alpha$, and only then when
% the residual $ F(\rho,\alpha) -\bar{f}_\alpha(-x+h-h_1) =0$. \newline
% \centerline{}

%\begin{proposition} 
%If $\hat{\rho}$ solves the minimization problem and 
%$J(\hat{\rho}) > 0,$ then at least one of the following conditions must be satisfied
%for almost every $x \in \Omega$:
%\begin{align*}
%\rho_r(x) =& \rho_{r_{0}} \\
%\rho_r(x) =& \rho_{r_{1}} \\
%\rho_{i}(x) =& \rho_{i_{0}}\\
%\rho_{i}(x) =& \rho_{i_{1}}. 
%\end{align*}
%\end{proposition}

\section{Numerical Results}
Approximate solutions of (\ref{minprob}) are sought through numerical
discretization and optimization.  The variational problem
(\ref{varprob}) was discretized with a first-order finite element
method, using piecewise bilinear elements on a uniform, rectangular
grid.  The design variable $\rho$ was approximated by a piecewise
constant function on the same uniform grid.  The nonlocal boundary
operators $T_\alpha$ defined by (\ref{Tdef}) were approximated
by explicitly calculating the Fourier coefficients of the traces
of the finite element basis, then truncating the sum in (\ref{Tdef}).
The resulting finite element scheme can be shown to converge and to 
conserve energy, provided
all the propagating terms are included in the sum \cite{Bao}.
This discretization leads to a large, sparse (except for the
boundary terms), non-Hermitian matrix problem, which for simplicity
is solved using the direct sparse solver in {\em Matlab}.

The integral in (\ref{minprob}) was approximated by a discrete
sum in $\alpha$.  By imposing $x$-axis symmetry in the designs, 
it suffices to integrate only over positive $\alpha$.  In the
following examples, we used 20 equally-spaced positive values
of $\alpha$ to approximate the integral.

Despite the convenience of imposing a positive lower bound on
the imaginary part of $\rho$ in Lemma~2.1 for obtaining a uniform
upper bound on solutions, we found that the numerical experiments
were quite insensitive to small dissipations. Thus in most of
the examples below, we set $\rho_{i_0} = 0$. 

After discretizing $J(\rho)$ through finite elements, optimization
was accomplished with a straightforward projected gradient descent 
algorithm as in \cite{Dobson}, using the adjoint as derived in 
Section~3 to calculate the gradient.
We performed a large number of numerical experiments with this
method, using different initial guesses for the design variable 
$\rho$, and varying the frequency $\omega$, source and focus 
locations $h$ and $h_1$,
and constraints $\rho_{r_0}, \rho_{r_1}, \rho_{i_0}, \rho_{i_1}$ 
on the real and imaginary parts of $\rho$.

Generally speaking, we found that the method was able to produce,
from almost any initial guess, a structure which produced
relatively high field intensity near the desired focus.  Some
parameter choices and initial guesses resulted in structures with
much better focusing properties than others.  

\begin{figure}
\centerline{
\begin{tabular}{c}
\includegraphics[width=60mm]{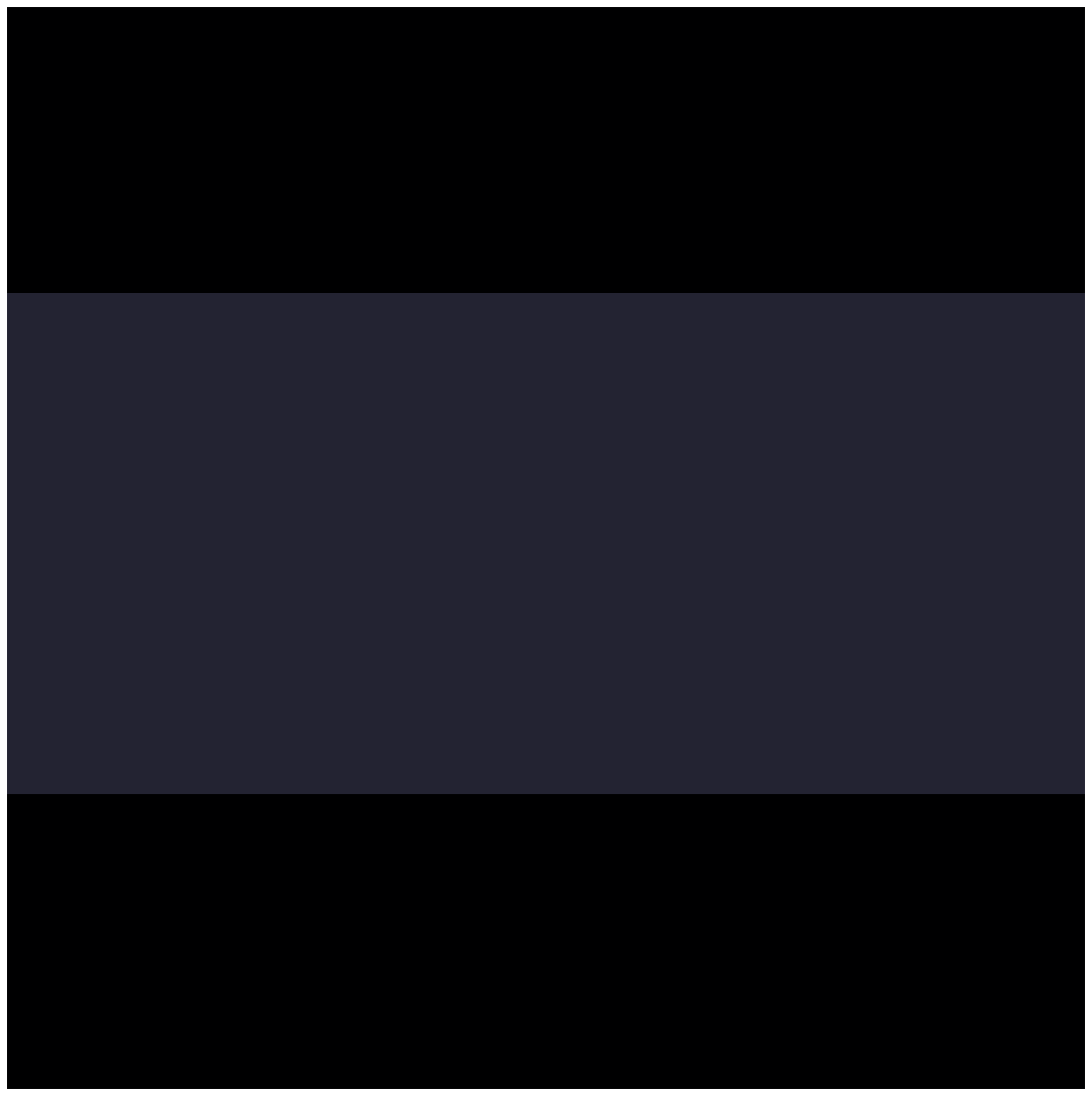}
\hspace*{4mm}
\includegraphics[width=60mm]{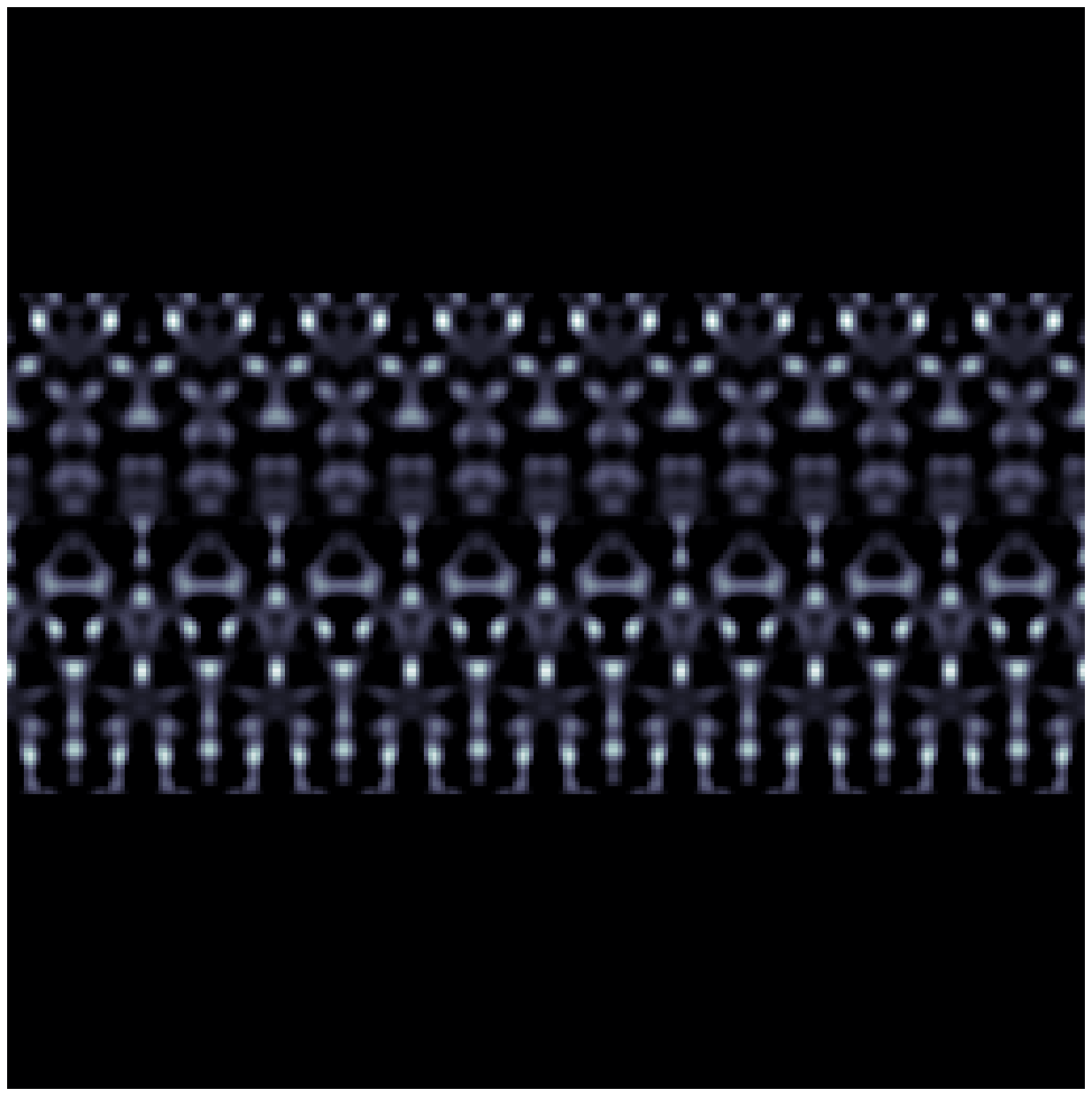} \\
\includegraphics[width=60mm]{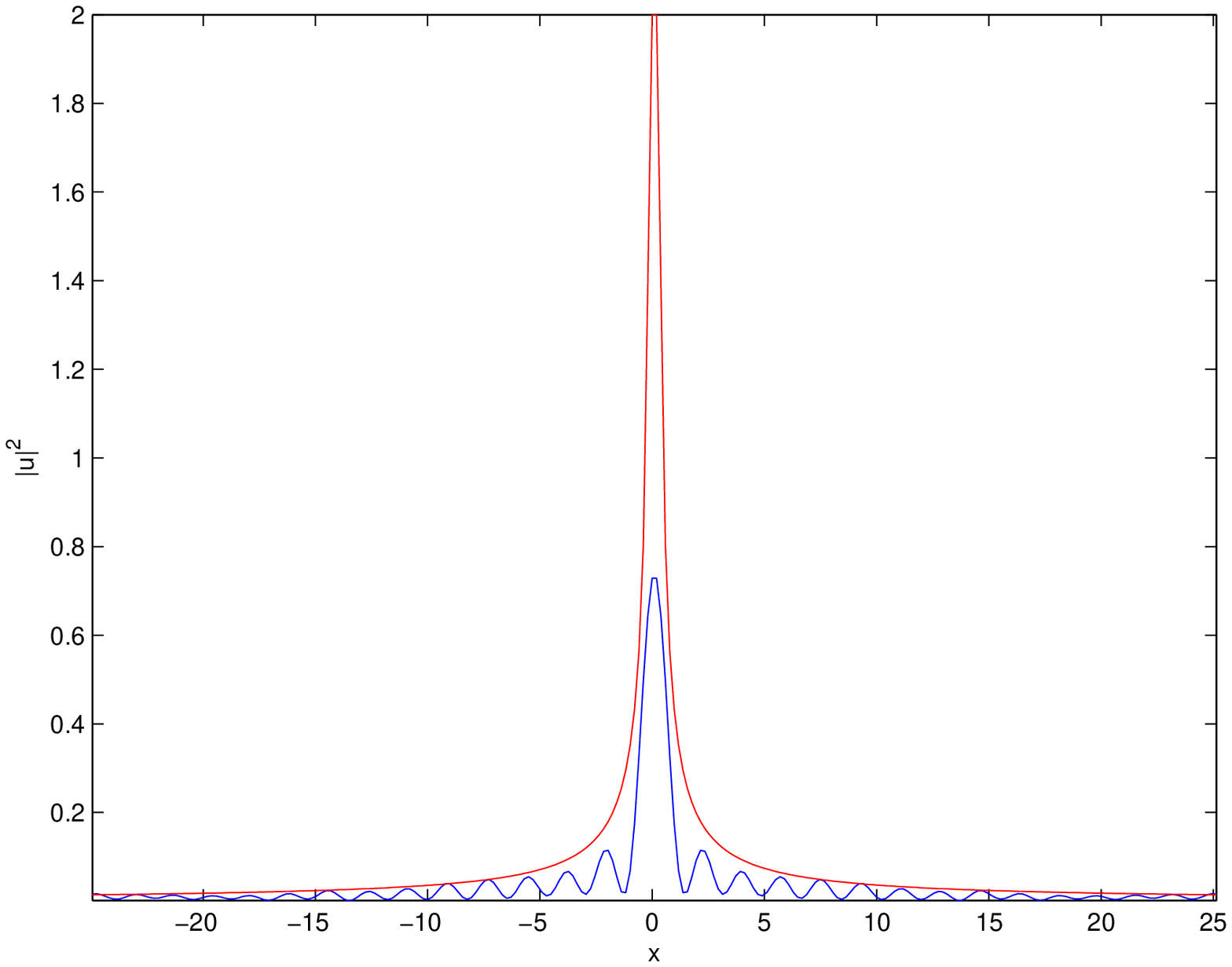} 
\hspace*{4mm}
\includegraphics[width=60mm]{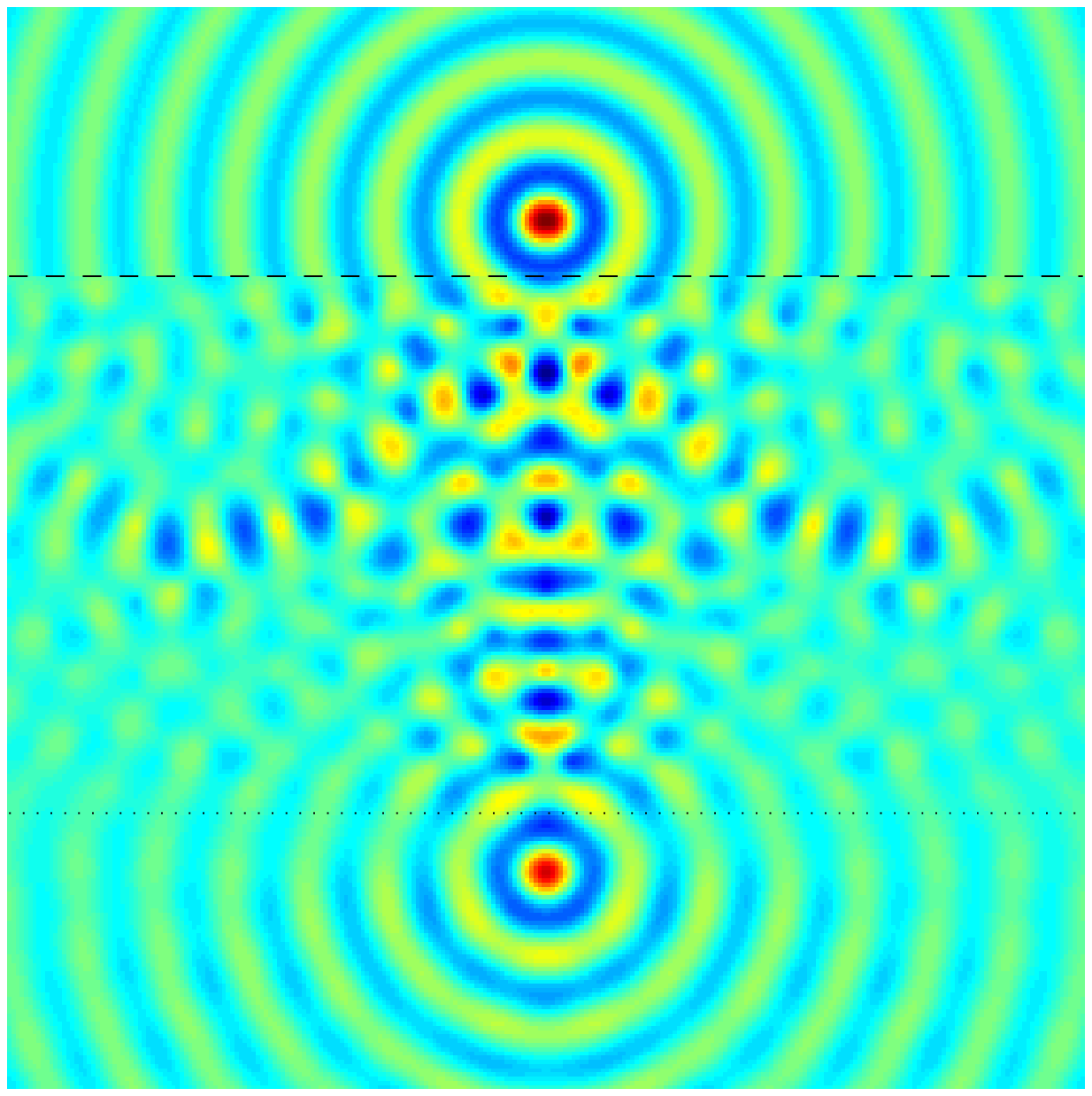}
\end{tabular}
}
\caption{Results from Example 1. Upper left: initial (real) $\rho$;
upper right: optimized solution; lower left: intensity cross section
(blue) versus point source (red).  Spot size is 0.424; lower right: 
real part of $E$ field
within the solution box.  Eight periods are shown.}
\label{example1}
\end{figure}

\begin{figure} \label{mode_plot}
\centerline{
\includegraphics[width=60mm]{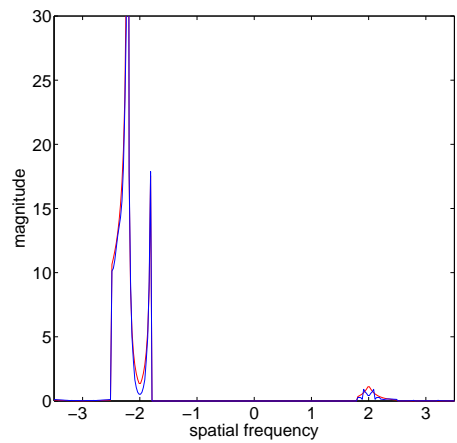}
}
\caption{Magnitude of evanescent modes for Example 1 (blue), versus
those of the target (red).}
\end{figure}

In the first experiment we start with a purely real material and 
allow the real part of the dielectric coefficient to vary between 
$\rho_{r_0} = 1$ and $\rho_{r_1}=12$. The source is positioned at 
$2.5$ units from the slab ($h=2.5$), and we are looking to obtain a focus 
$2.5$ units on the opposite side of the slab ($h_1=2.5$). Through numerical optimization
we discover the structure shown in Figure \ref{example1} which gives a spot size 
$0.424 \lambda$. 
True subwavelength imaging is only possible if evanescent modes are present
at the interface between the structure and the transmission medium.  
Figure~\ref{mode_plot} shows the evanescent modes for the generated image versus those of the
target, clearly showing that the structure produces evanescent modes 
which approximate those of the objective.

\begin{figure}
\centerline{
\begin{tabular}{c}
\includegraphics[width=60mm]{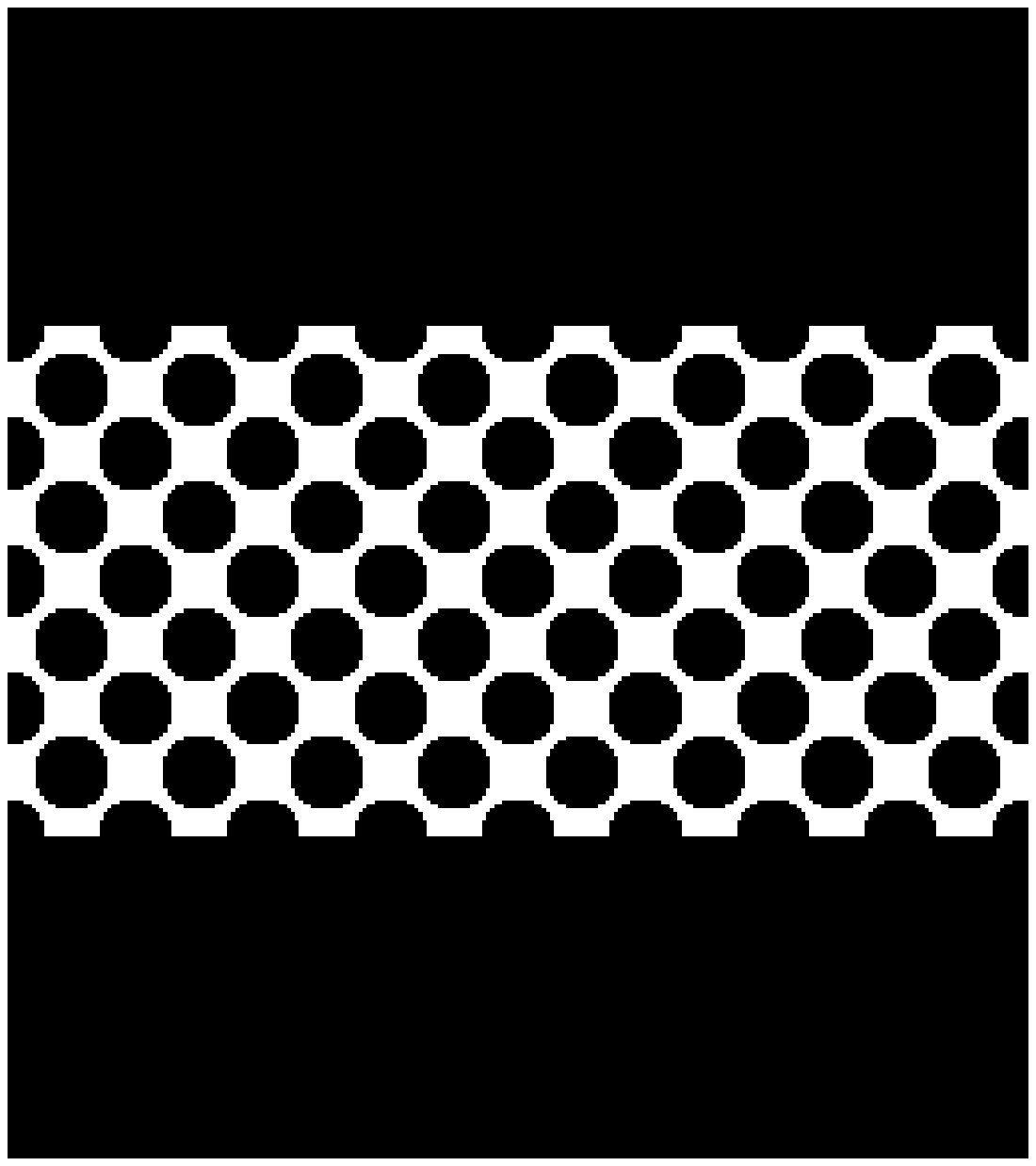}
\hspace*{4mm}
\includegraphics[width=60mm]{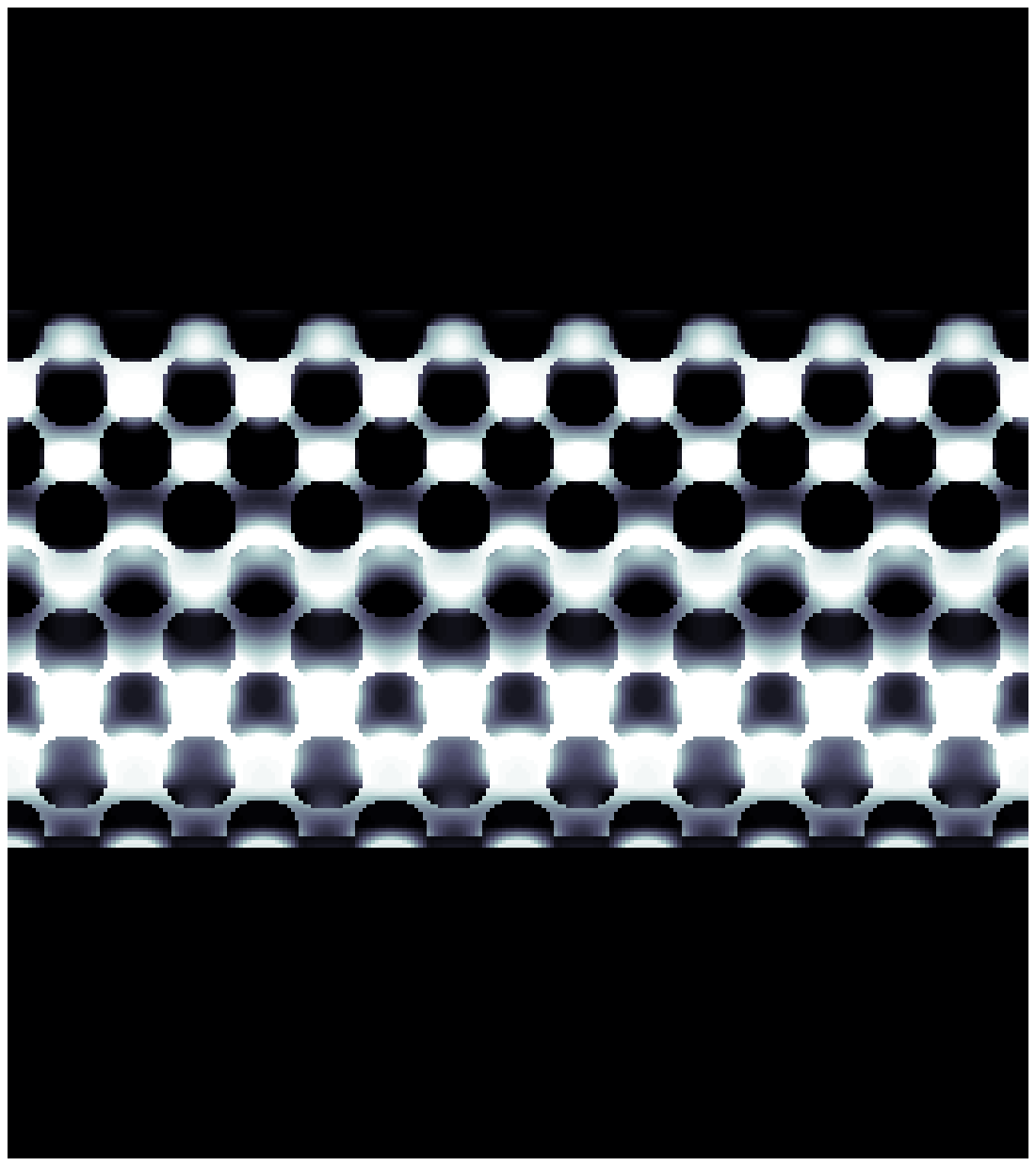} \\
\includegraphics[width=60mm]{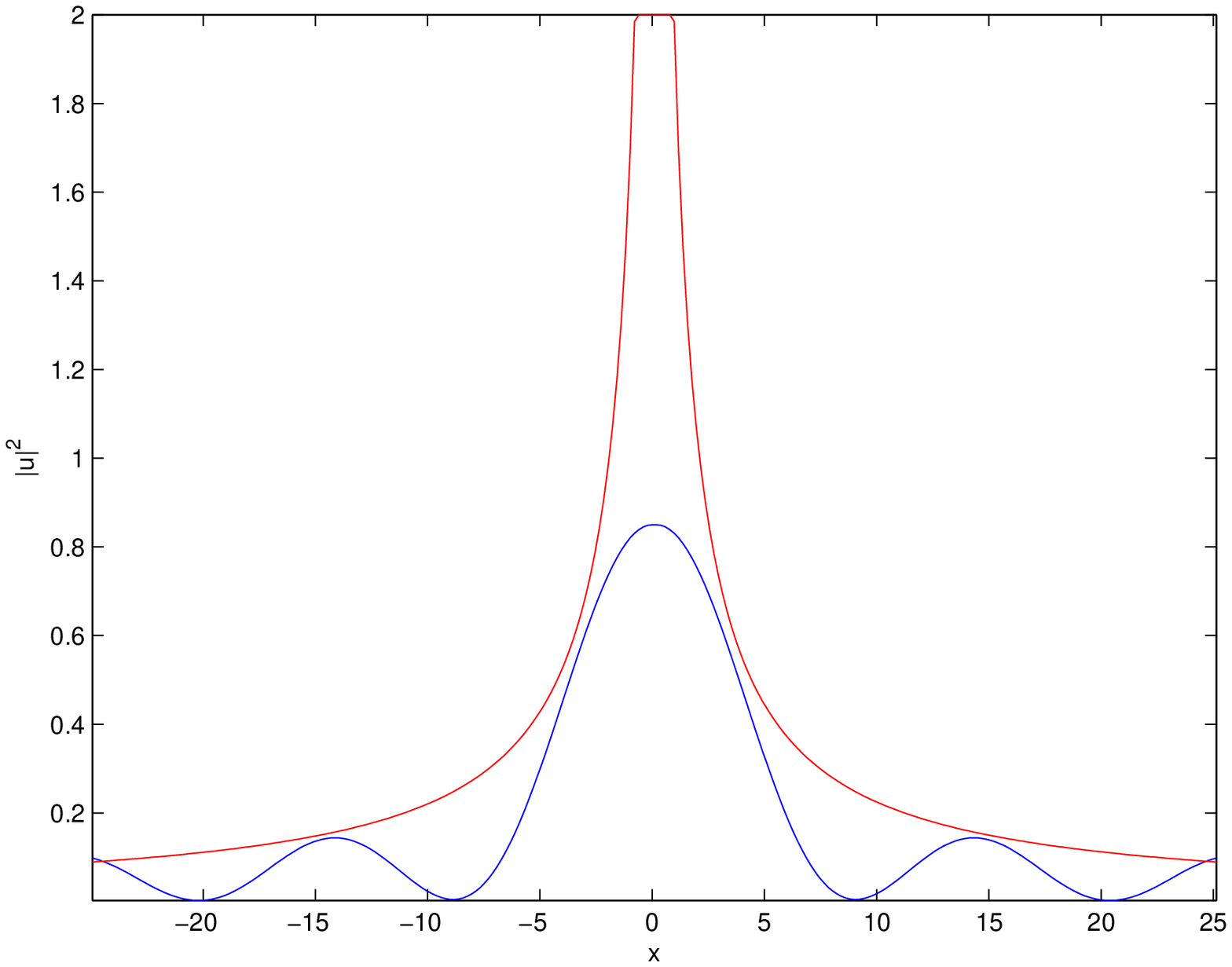} 
\hspace*{4mm}
\includegraphics[width=60mm]{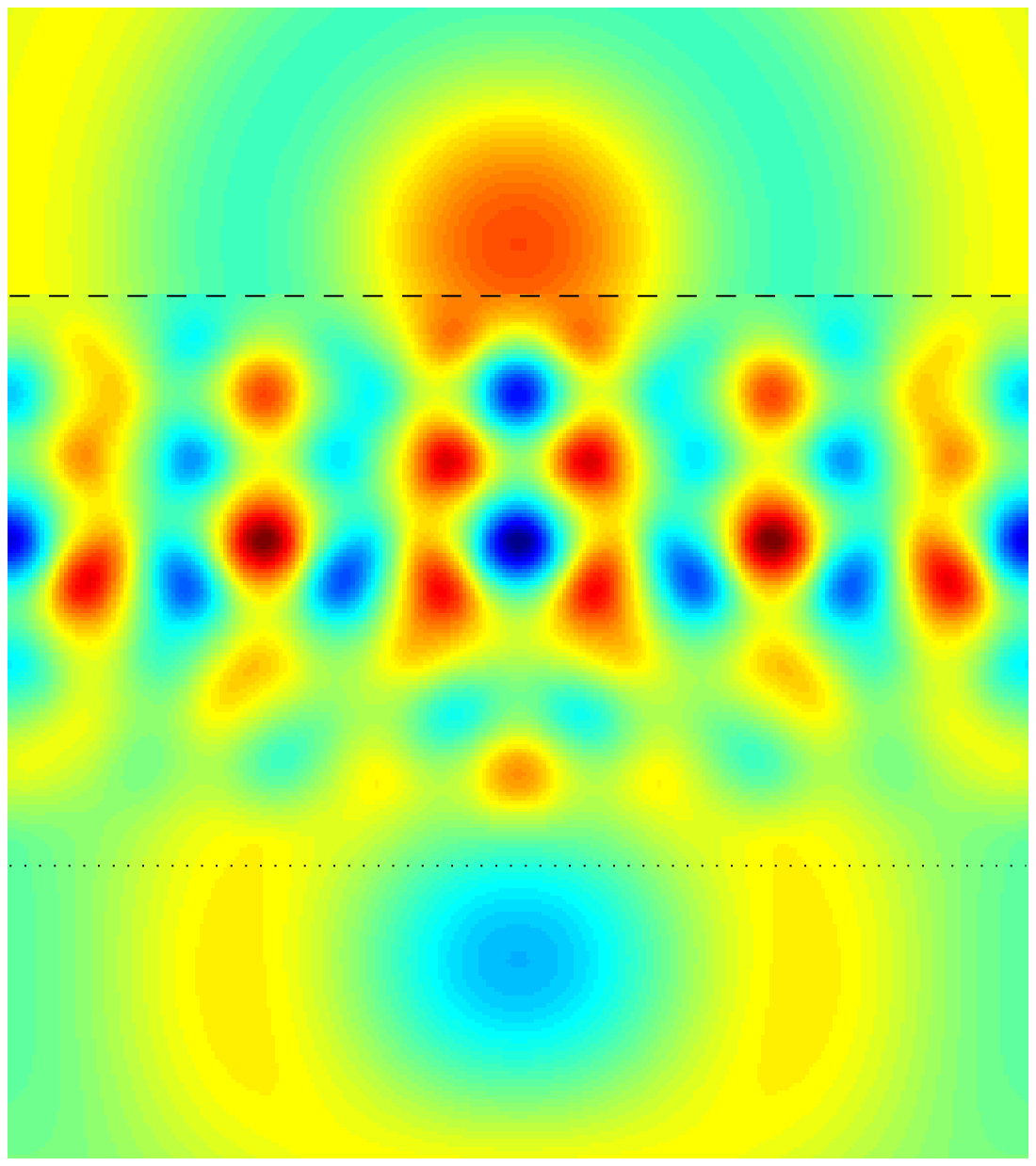}
\end{tabular}
}
\caption{Results from Example 3. Upper left: initial (real) $\rho$;
upper right: optimized solution; lower left: intensity cross section
(blue) versus point source (red).  Spot size is 0.395; lower right: 
real part of $E$ field
within the solution box.  Eight periods are shown.}
\label{example3}
\end{figure}

In the second example we image a source far away from the lens ($h=90$), and we are looking for 
the lens that will produce the best image ten units from the lens ($h_1=10$). Far away objects
are much more difficult to image, but through optimization of the structure, we are able to 
obtain a spot of size $0.38 \lambda$. 
In the third example we start with a photonic crystal and allow purely real structures 
varying between $\rho_{r_0} = 1$ and $\rho_{r_1}=12$. The source is positioned at 
$2.4$ units from the slab ($h=2.4$) and we are looking to obtain a focus 
four units away on the opposite side of the slab ($h_1=4$). The optimized structure is shown in Figure \ref{example3},
and it gives a focus with spot size $0.395 \lambda$, which is much 
better than the one obtained if we just use photonic crystal as suggested by Luo \textit{et al.} 
which produced a focus with spot size $0.67 \lambda$ \cite{Luo}.

\begin{figure}
\centerline{
\begin{tabular}{c}
\includegraphics[width=60mm]{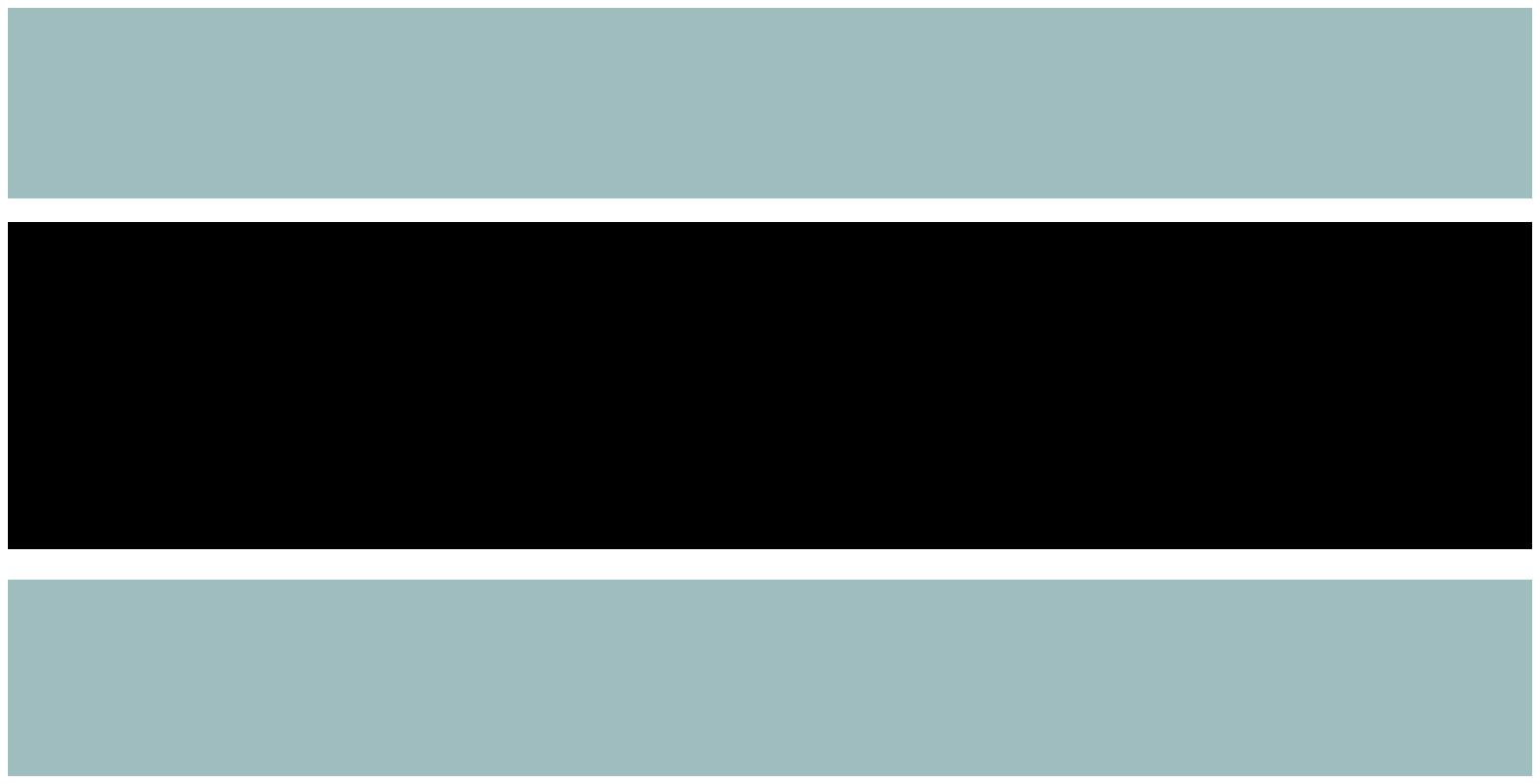}
\hspace*{4mm}
\includegraphics[width=60mm]{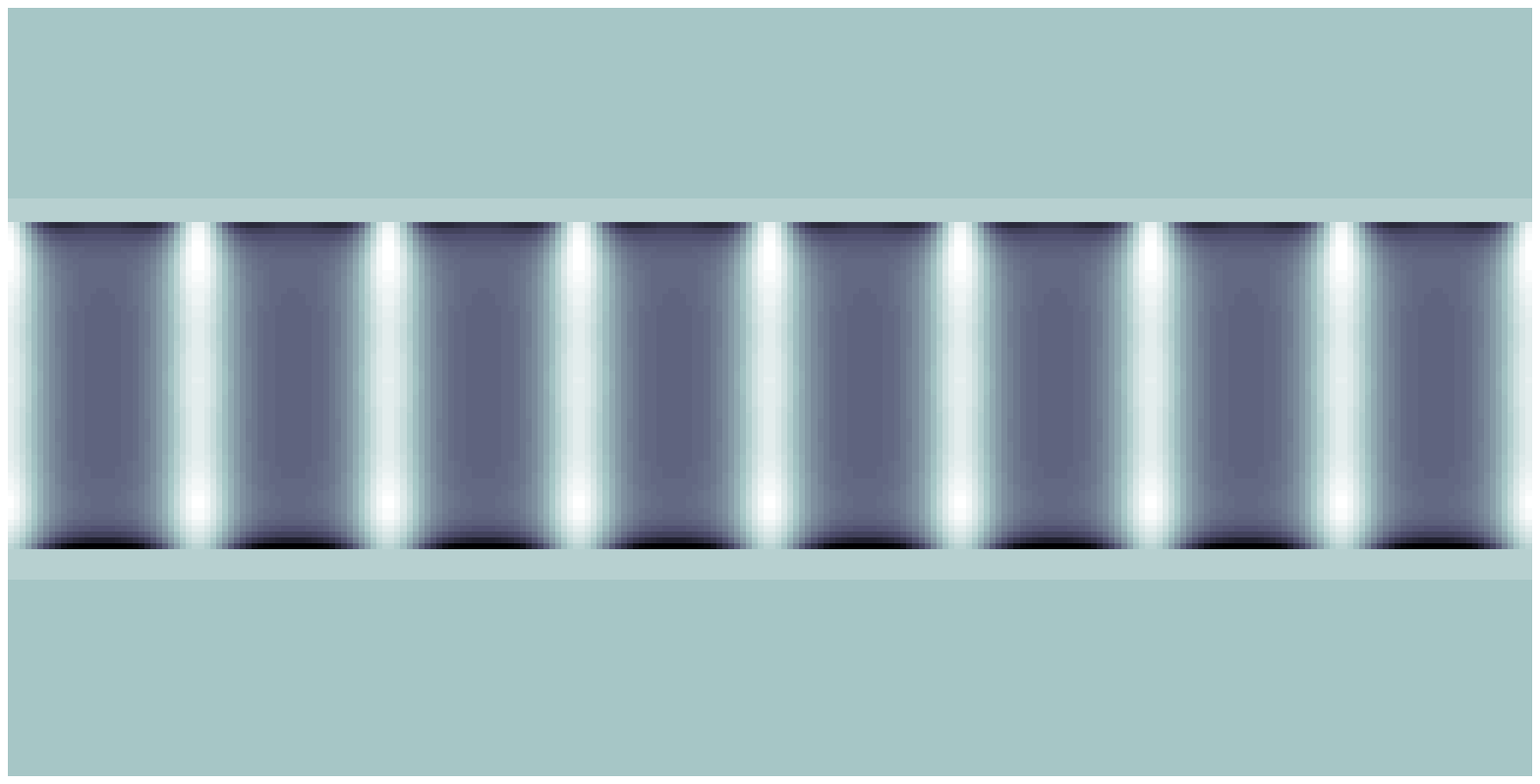} \\
\includegraphics[width=60mm]{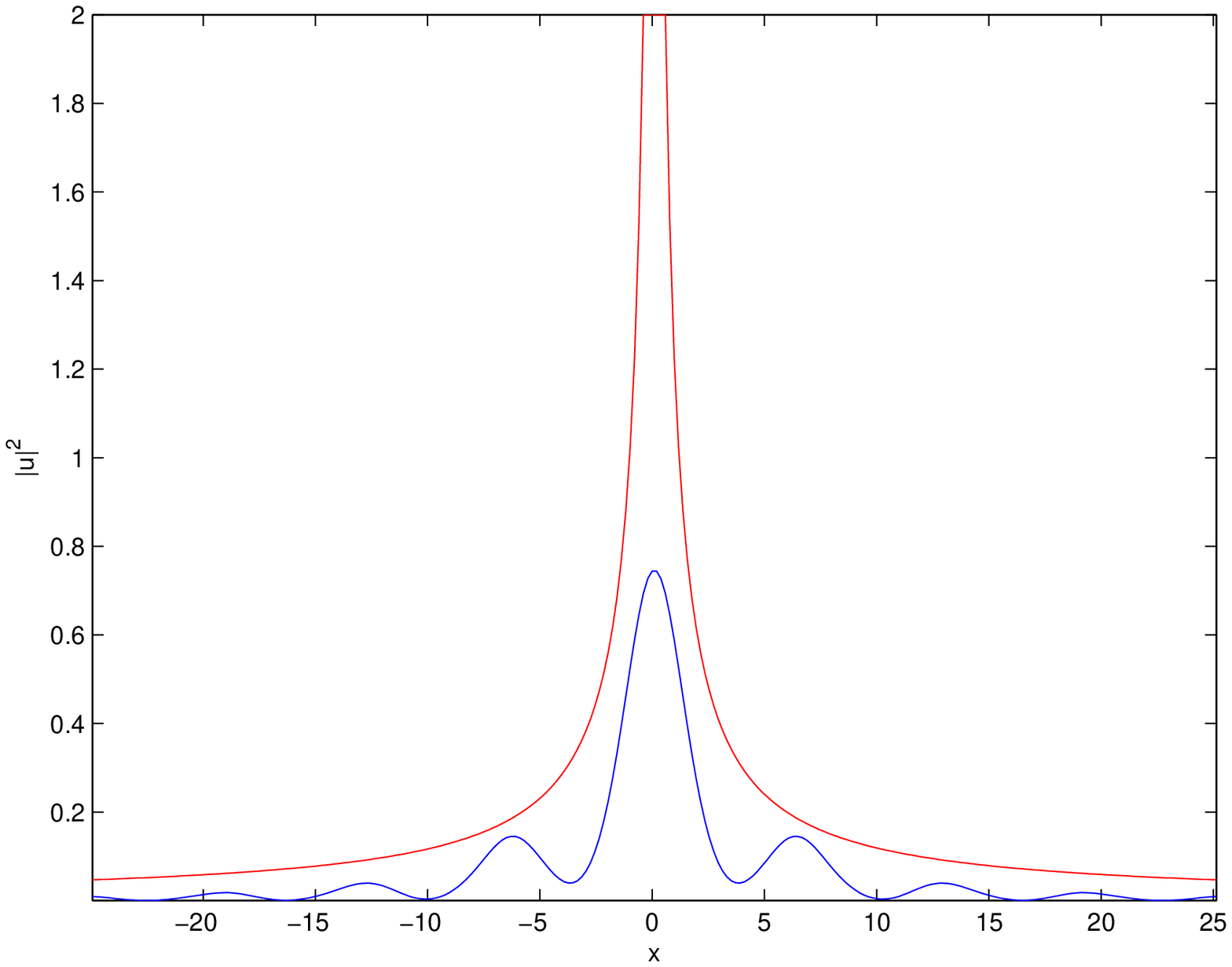} 
\hspace*{4mm}
\includegraphics[width=60mm]{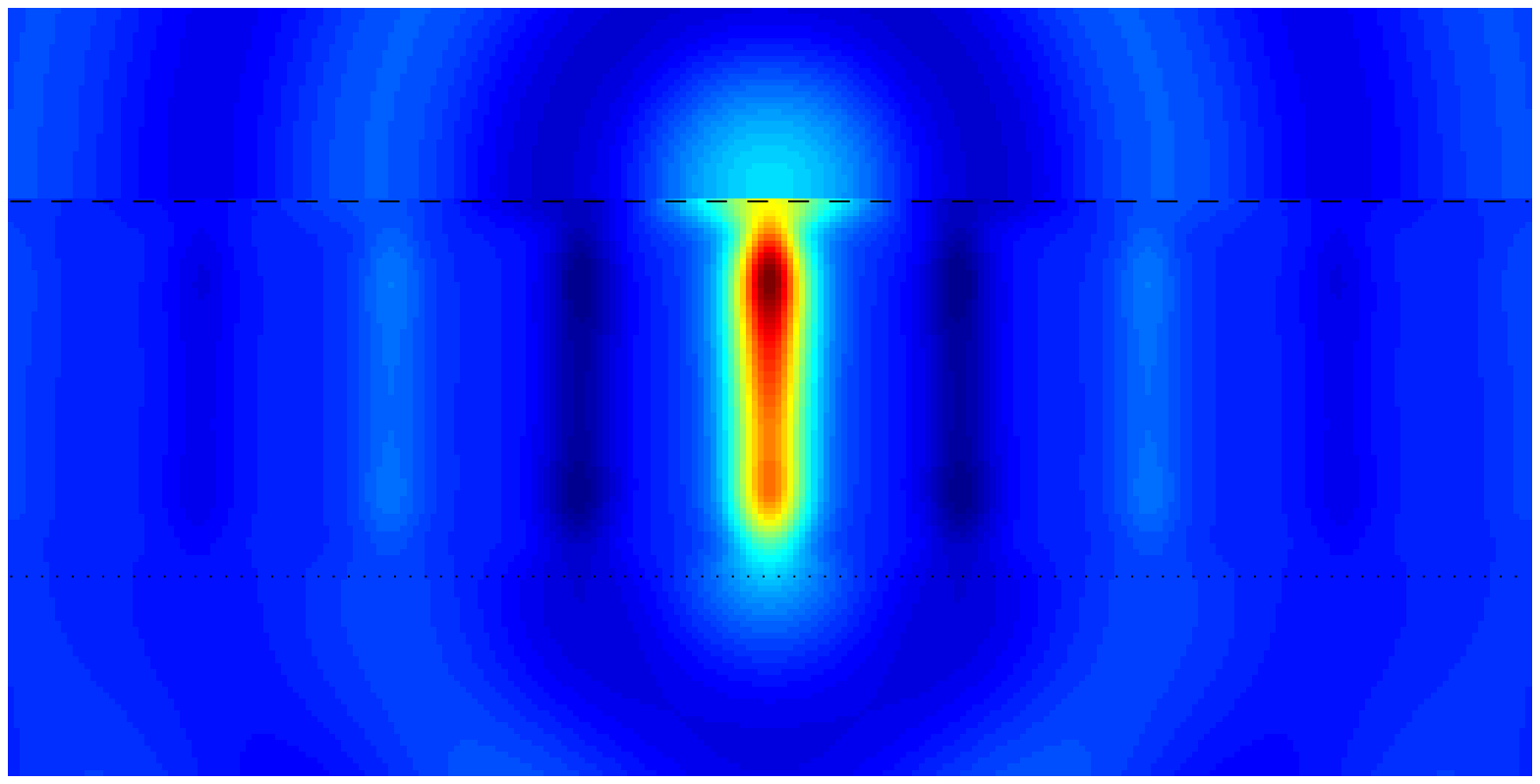}
\end{tabular}
}
\caption{Results from Example 4. Upper left: initial (real) $\rho$;
upper right: optimized solution; lower left: intensity cross section
(blue) versus point source (red).  Spot size is 0.284; lower right: 
real part of $E$ field
within the solution box.  Eight periods are shown.}
\label{example4}
\end{figure}

Example four allows both the real and imaginary parts of the dielectric coefficient to vary 
($\rho_{r_0} = 1$, $\rho_{r_1}=12$, and $\rho_{i_0} = 0$, $\rho_{i_1}=1$). 
The distance between the source and the lens and
the distance between the lens and the image are set to four units ($h=4$, $h_1=4$). 
The optimized structure gives a focus with a 
spot size $0.284 \lambda$ (Figure \ref{example4}), which is a significant improvement to those obtained by 
structures described in the research literature so far, although as far as
we know, no real materials exist with these dielectric coefficients.
%In the fifth example we decrease the distance between the source and the lens and
%the distance between the lens and the image to $0.9$ units ($h=0.9$, $h_1=0.9$). 
%We start with a purely real material and 
%allow the real part of the dielectric coefficient to vary between 
%$\rho_{r_0} = 1$ and $\rho_{r_1}=12.5$. Through optimization we obtain a focus with
%spot size $0.397 \lambda$ (Figure \ref{example5}).
 
All of the numerical
experiments were computationally intensive.  Each iteration
required the solution of a family of diffraction problems, two for each
$\alpha$, and because of the crude optimization method employed,
many iterations were typically required.  Most of the examples
below took on the order of two days to run on a workstation.
Since our purpose here was simply to illustrate the feasibility
of designing such structures through mathematical optimization,
we did not devote much effort to improving the efficiency of the
numerical methods.  We believe that computation time could be
improved by at least an order of magnitude with existing, but
more sophisticated, numerical methods.

\begin{figure}
\centerline{
\begin{tabular}{c}
\includegraphics[width=60mm]{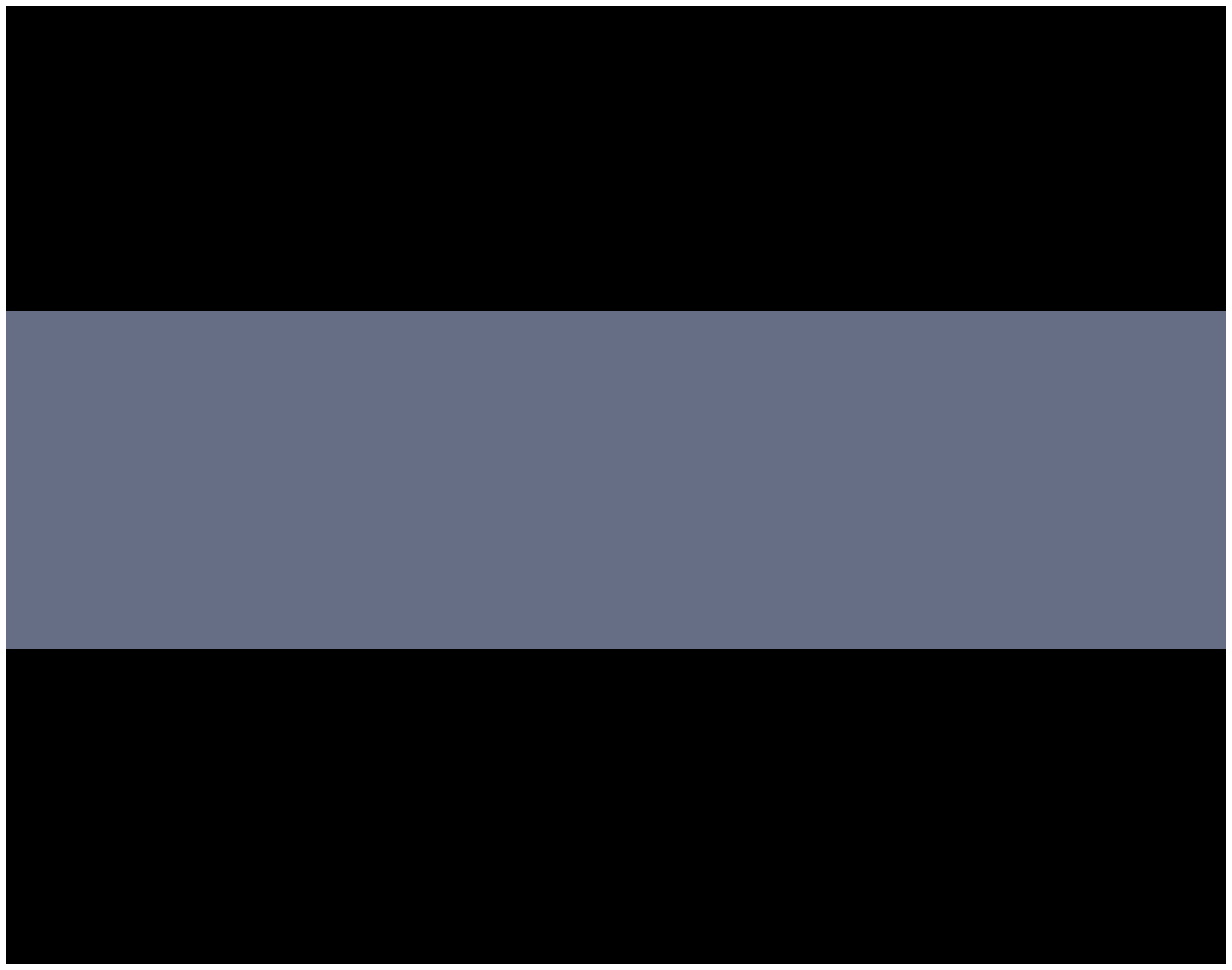}
\hspace*{4mm}
\includegraphics[width=60mm]{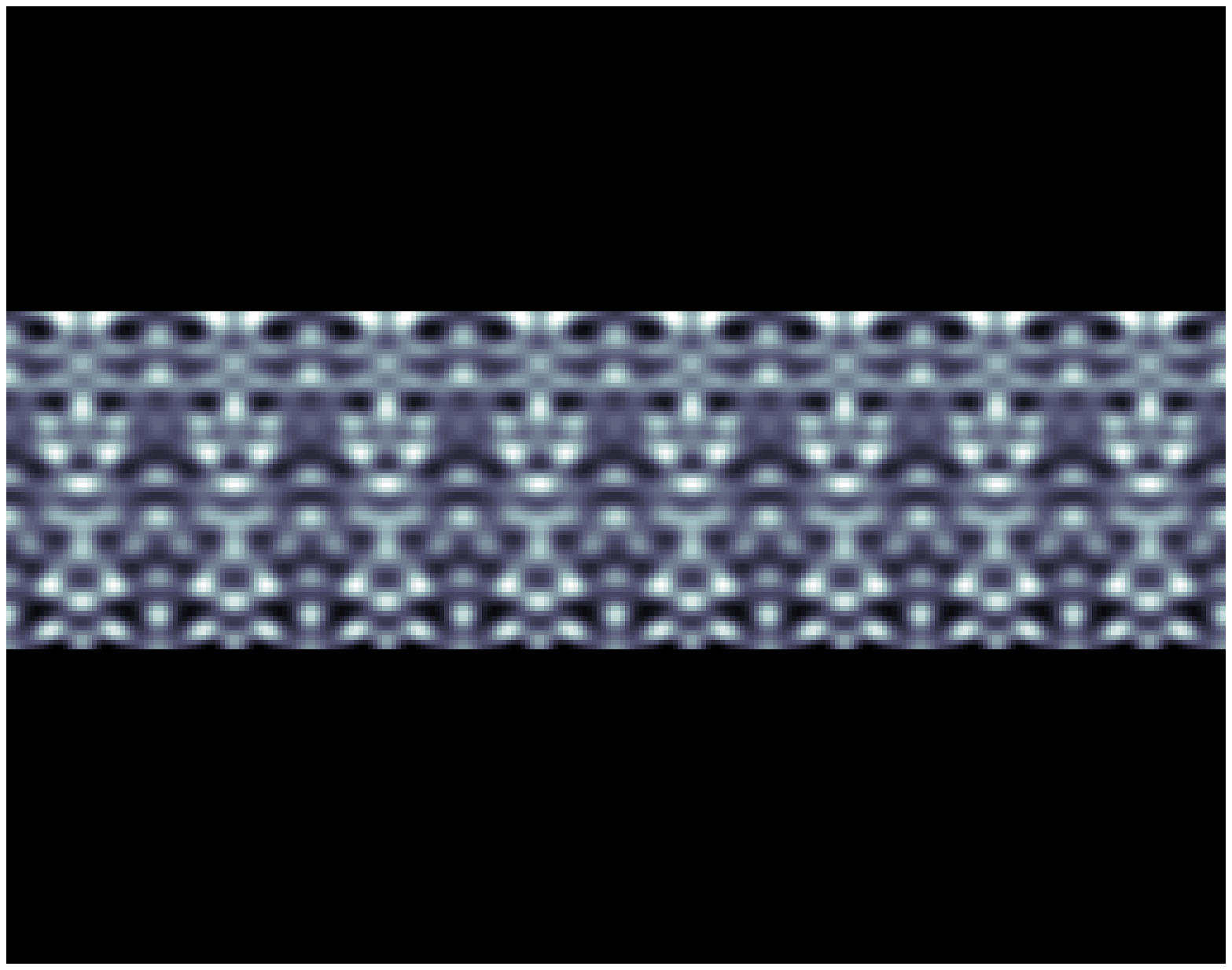} \\
\includegraphics[width=60mm]{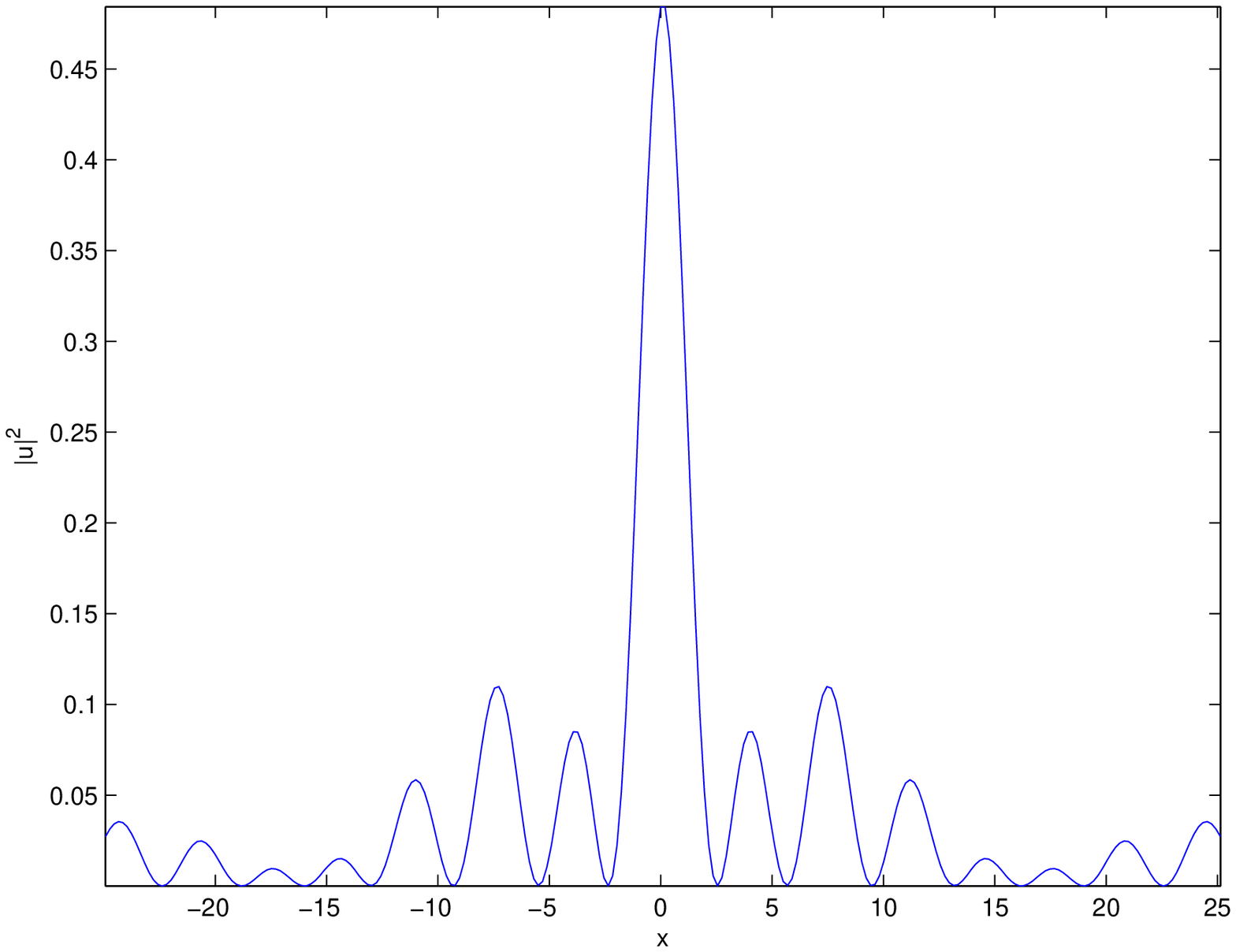} 
\hspace*{4mm}
\includegraphics[width=60mm]{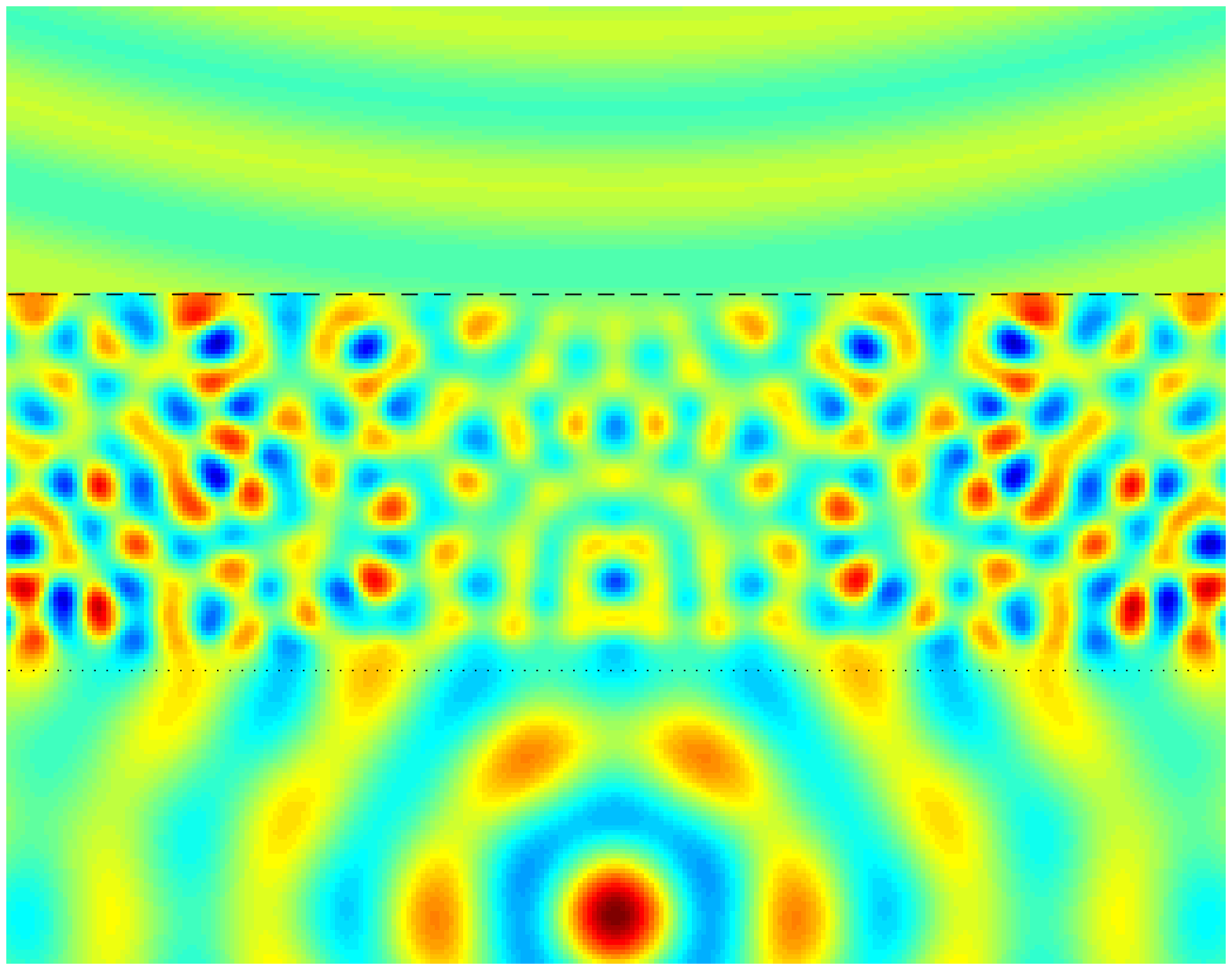}
\end{tabular}
}
\caption{Results from Example 2. Upper left: initial (real) $\rho$;
upper right: optimized solution; lower left: intensity cross section.  
Spot size is 0.380; lower right: real part of $E$ field
within the solution box.  Eight periods are shown.}
\label{example2}
\end{figure}

%\begin{figure}
%\centerline{
%\begin{tabular}{c}
%\includegraphics[width=60mm]{ex5_eps_initial.eps}
%\hspace*{4mm}
%\includegraphics[width=60mm]{ex5_epsilon.eps} \\
%\includegraphics[width=60mm]{ex5_spot.eps} 
%\hspace*{4mm}
%\includegraphics[width=60mm]{ex5_realu.eps}
%\end{tabular}
%}
%\caption{Results from Example 5. Upper left: initial (real) $\rho$;
%upper right: optimized solution; lower left: intensity cross section
%(blue) versus point source (red).  Spot size is 0.397; lower right: 
%real part of $E$ field
%within the solution box.  Eight periods are shown.}
%\label{example5}
%\end{figure}

\section{Conclusions}
We have demonstrated the feasibility of designing periodic structures
with subwavelength focusing properties, via mathematical optimization.
The approach is mathematically sound and through numerical 
discretization yields plausible, if somewhat non-intuitive, solutions.

We point out two weaknesses with the approach presented here, all
of which we believe could be improved through further work.
First, the performance of the
optimized structures tends to be very sensitive to small perturbations
in material parameters.  This fact combined with the relative
complexity of the solutions means that attempting to fabricate 
these structures is at this point not an attractive idea.
Both the sensitivity and the complexity could be addressed through
adding appropriate constraints or penalties to the objective. For
example, optimization through a level-set approach would
necessarily yield structures composed of only two materials,
with no intermediate-index areas \cite{Santosa}.
This together with total variation penalties may
yield ``simpler'' solutions.  Of course such constraints
may incur a decrease in the performance of the structure.

Second, the designs created with this approach are not translationally
invariant.  Specifically, moving the structure laterally relative
to the point source may result in decrease of focus.  There are a few
ways to circumvent this problem.  In the simplest case, the distance
$h$ from the point source to the structure is large so that the
wavefront impinging on the structure is nearly planar.  Our experiments
have shown that in this case focusing is not dependent on the lateral
position of the point source (and is much less senstive to vertical
translations), although the focus does translate periodically with
the structure.  One can also optimize for structures whose
period is very small compared to the wavelength, although such
structures tend to have less focusing power.  The ultimate solution
to this problem would require building the translation invariance 
of the focus (not of the medium) into the objective function.
Unfortunately this would increase the complexity of the computations
beyond the capability of the simplistic approach presented here,
but with further work and refinement it should be possible.

\end{document}